\documentclass{amsart}
\usepackage{amsmath,amssymb,amsfonts}
\usepackage{eucal}
\usepackage{verbatim}

\usepackage{enumerate}
\newenvironment{enumroman}{\begin{enumerate}[\upshape (i)]}
                                                {\end{enumerate}}

\numberwithin{equation}{section}

\theoremstyle{plain}
\newtheorem{theorem}[equation]{Theorem}

\newtheorem{cor}[equation]{Corollary}
\newtheorem{prop}[equation]{Proposition}
\newtheorem{lemma}[equation]{Lemma}

\theoremstyle{definition}
\newtheorem{definition}[equation]{Definition}

\newtheorem*{note}{Note}

\newtheorem*{thank}{Acknowledgments}

\input xy
\xyoption{all}

\newcommand{\Deltaop}{{\bf \Delta}^{op}}
\newcommand{\colim}{\text{colim}}

\newcommand{\nerve}{\text{nerve}}
\newcommand{\Hom}{\text{Hom}}
\newcommand{\SSets}{\mathcal{SS}ets}

\newcommand{\SSetsd}{\mathcal{SS}ets^{\mathcal D}}

\newcommand{\LSSets}{\mathcal{LSS}ets}

\newcommand{\Secat}{\mathcal Se \mathcal Cat}
\newcommand{\Ho}{\text{Ho}}
\newcommand{\map}{\text{map}}
\newcommand{\Map}{\text{Map}}
\newcommand{\ob}{\text{ob}}
\newcommand{\hoequiv}{\text{hoequiv}}

\newcommand{\sco}{\mathcal{SC}_\mathcal O}
\newcommand{\css}{\mathcal{CSS}}
\newcommand{\cosk}{\text{cosk}}
\newcommand{\sk}{\text{sk}}

\newcommand{\sesp}{\mathcal Se \mathcal Sp}
\newcommand{\xu}{{\underline x}}

\newcommand{\yu}{{\underline y}}

\begin{document}

\title{Three Models for the Homotopy Theory of Homotopy Theories}

\author[J.E. Bergner]{Julia E. Bergner}

\address{Department of Mathematics \\
University of Notre Dame \\
Notre Dame, IN 46556}

\curraddr{Kansas State University \\
138 Cardwell Hall \\
Manhattan, KS 66506}

\email{bergnerj@member.ams.org}

\date{\today}



\thanks{The author was partially supported by a Clare Boothe Luce
Foundation Graduate Fellowship.}

\begin{abstract}
Given any model category, or more generally any category with weak
equivalences, its simplicial localization is a simplicial category
which can rightfully be called the ``homotopy theory" of the model
category. There is a model category structure on the category of
simplicial categories, so taking its simplicial localization
yields a ``homotopy theory of homotopy theories." In this paper we
show that there are two different categories of diagrams of
simplicial sets, each equipped with an appropriate definition of
weak equivalence, such that the resulting homotopy theories are
each equivalent to the homotopy theory arising from the model
category structure on simplicial categories.  Thus, any of these
three categories with the respective weak equivalences could be
considered a model for the homotopy theory of homotopy theories.
One of them in particular, Rezk's complete Segal space model
category structure on the category of simplicial spaces, is much
more convenient from the perspective of making calculations and
therefore obtaining information about a given homotopy theory.
\end{abstract}

\maketitle

\section{Introduction}
Classical homotopy theory considers topological spaces, up to weak
homotopy equivalence. Eventually, the structure of the category of
topological spaces making it possible to talk about its ``homotopy
theory" was axiomatized; it is known as a model category
structure.  In particular, given a model category structure on an
arbitrary category, we can talk about its homotopy category.  More
generally, we can think about the ``homotopy theory" given by that
category with its particular class of weak equivalences, where the
homotopy theory encompasses the homotopy category as well as
higher-order information. One might ask what specifically is meant
by a homotopy theory.

One answer to this question uses simplicial categories, which in
this paper we will always take to mean categories enriched over
simplicial sets.  Given a model category $\mathcal M$, taking its
simplicial localization with respect to its subcategory of weak
equivalences yields a simplicial category $\mathcal{LM}$
\cite[4.1]{dk}. The simplicial localization encodes the known
homotopy-theoretic information of the model category, so one point
of view is that this simplicial category is the homotopy theory
associated to the model category structure.  Set-theoretic issues
aside, we can also construct the simplicial localization for any
category with a subcategory of weak equivalences, so therefore we
can speak of an associated homotopy theory even in this more
general situation.

Given two homotopy theories, one can ask whether they are
equivalent to one another in some natural sense. There is a notion
of weak equivalence between two simplicial categories which is a
simplicial analogue of an equivalence between categories. These
weak equivalences are known as \emph{DK-equivalences}, where the
``DK" refers to the fact that they were first defined by Dwyer and
Kan in \cite{dk1}. In fact, there is a model category structure
$\mathcal{SC}$ on the category of all (small) simplicial
categories in which the weak equivalences are these
DK-equivalences \cite[1.1]{simpcat}. The associated homotopy
theory of simplicial categories is what we will refer to as the
homotopy theory of homotopy theories.

In \cite{rezk}, Rezk takes steps toward finding a model other than
that of simplicial categories for the homotopy theory of homotopy
theories.  He defines complete Segal spaces, which are simplicial
spaces satisfying some nice properties (Definitions
\ref{SegalSpace} and \ref{CSegalSpace} below) and constructs a
functor which assigns a complete Segal space to any simplicial
category. He considers a model category structure $\mathcal{CSS}$
on the category of all simplicial spaces in which the weak
equivalences are levelwise weak equivalences of simplicial sets
and then localizes it in such a way that the local objects are the
complete Segal spaces (Theorem \ref{CSS}).

However, Rezk does not construct a functor from the category of
complete Segal spaces to the category of simplicial categories,
nor does he discuss the model category $\mathcal{SC}$. In this
paper, we complete his work by showing that $\mathcal{SC}$ and
$\mathcal{CSS}$ have equivalent homotopy theories. This result is
helpful in that the weak equivalences between complete Segal
spaces are easy to identify (see Proposition \ref{DKReedy} below),
unlike the weak equivalences between simplicial categories, and
therefore making any kind of calculations would be much easier in
$\mathcal{CSS}$. Using terminology of Dugger \cite{dug}, this
model category $\mathcal{CSS}$ is a \emph{presentation} for the
homotopy theory of homotopy theories, since it is a localization
of a category of diagrams of spaces.

In order to prove this result, we make use of an intermediate
category.  Consider the full subcategory $\Secat$ of the category
of simplicial spaces whose objects are simplicial spaces with a
discrete simplicial set in degree zero.  We will prove the
existence of two model category structures on $\Secat$, each with
the same class of weak equivalences.  The first of these
structures, which we denote $\Secat_c$, has as cofibrations the
maps which are levelwise cofibrations of simplicial sets.  (An
alternate proof of the existence of this model category structure
is given by Hirschowitz and Simpson \cite[2.3]{hs}.  They actually
prove the existence of such a model category structure for Segal
$n$-categories, whereas we consider only the case where $n=1$.)
The second model category structure, which we denote $\Secat_f$,
has as fibrations maps which can be thought of as localizations of
levelwise fibrations of simplicial sets, although strictly
speaking they cannot be obtained in this way. We use these model
category structures to produce a chain of Quillen equivalences
\[ \mathcal{SC} \leftrightarrows \Secat_f \rightleftarrows \Secat_c
\rightleftarrows \mathcal{CSS}. \]  (In each case, the topmost
arrow is the left adjoint of the adjoint pair.) Notice that we can
obtain a single Quillen equivalence $\Secat_f \rightleftarrows
\mathcal {CSS}$ via composition.  Since Quillen equivalent model
categories have DK-equivalent simplicial localizations
(Proposition \ref{qedk}), all three of these categories with their
respective weak equivalences give models for the homotopy theory
of homotopy theories.

\subsection{Organization of the Paper}
We begin in section 2 by recalling standard information about
model category structures and simplicial objects.  In section 3,
we state the definitions of simplicial categories, complete Segal
spaces, and Segal categories, and we give some basic results about
each. In section 4, we set up some constructions on Segal
precategories that we will need in order to prove our model
category structures. In section 5, we prove the existence of a
model category structure $\Secat_c$ on the category of Segal
precategories which we then in section 6 show is Quillen
equivalent to Rezk's complete Segal space model category structure
$\css$. In section 7, we prove the existence of the model category
structure $\Secat_f$ on the category of Segal precategories and
prove that it is Quillen equivalent to $\Secat_c$.  We then show
in section 8 that $\Secat_f$ is Quillen equivalent to the model
category structure $\mathcal{SC}$ on simplicial categories.
Section 9 contains the proofs of some technical lemmas.

\begin{thank}
This paper is a version of my Ph.D. thesis at the University of
Notre Dame \cite{thesis}.  I would like to thank my advisor Bill
Dwyer for his help and encouragement on this paper.  I would also
like to thank Charles Rezk, Bertrand To\"{e}n, and the referee for
helpful comments.
\end{thank}

\section{Background on Model Categories and Simplicial Objects}

\subsection{Model Categories}
Recall that a model category structure on a category $\mathcal C$
is a choice of three distinguished classes of morphisms:
fibrations ($\twoheadrightarrow$), cofibrations
($\hookrightarrow$), and weak equivalences ($\tilde \rightarrow$).
A (co)fibration which is also a weak equivalence is an
\emph{acyclic (co)fibration}. With this choice of three classes of
morphisms, $\mathcal C$ is required to satisfy five axioms MC1-MC5
which can be found in \cite[3.3]{ds}.

In all the model categories we use, the factorizations given by
axiom MC5 can be chosen to be functorial \cite[1.1.1]{hovey}. An
object $X$ in a model category is \emph{fibrant} if the unique map
$X \rightarrow \ast$ to the terminal object is a fibration.
Dually, $X$ is \emph{cofibrant} if the unique map from the initial
object $\phi \rightarrow X$ is a cofibration.  Given any object
$X$, the functorial factorization of  the map $X \rightarrow \ast$
as the composite of an acyclic cofibration followed by a fibration
\[ \xymatrix@1{X \text{ } \ar@{^{(}->}[r]^\sim & X^f \ar@{->>}[r] & \ast}
\] gives us the object $X^f$, the \emph{fibrant replacement} of
$X$.  Dually, we can define its \emph{cofibrant replacement} $X^c$
using the functorial factorization
\[ \xymatrix@1{\phi \text{ } \ar@{^{(}->}[r] & X^c \ar@{->>}[r]^\sim & X}. \]

All the model category structures that we work with are
cofibrantly generated.  In a cofibrantly generated model category,
there are two sets of specified morphisms, the generating
cofibrations and the generating acyclic cofibrations, such that a
map is an acyclic fibration if and only if it has the right
lifting property with respect to the generating cofibrations, and
a map is a fibration if and only if it has the right lifting
property with respect to the generating acyclic cofibrations
\cite[11.1.2]{hirsch}.  To prove that a particular category with a
choice of weak equivalences has a cofibrantly generated model
category structure, we need the following definition.

\setcounter{equation}{1}

\begin{definition} \cite[10.5.2]{hirsch}
Let $\mathcal C$ be a category and $I$ a set of maps in $\mathcal
C$.  Then an $I$-\emph{injective} is a map which was the right
lifting property with respect to every map in $I$.  An
$I$-\emph{cofibration} is a map with the left lifting property
with respect to every $I$-injective.
\end{definition}

We are now able to state the theorem that we use in this paper to
prove the existence of specific model category structures.

\begin{theorem} \cite[11.3.1]{hirsch} \label{CofGen}
Let $\mathcal M$ be a category with a specified class of weak
equivalences which satisfies model category axioms MC1 and MC2.
Suppose further that the class of weak equivalences is closed
under retracts.  Let $I$ and $J$ be sets of maps in $\mathcal M$
which satisfy the following conditions:
\begin{enumerate}
\item Both $I$ and $J$ permit the small object argument
\cite[10.5.15]{hirsch}.

\item Every $J$-cofibration is an $I$-cofibration and a weak
equivalence.

\item Every $I$-injective is a $J$-injective and a weak
equivalence.

\item One of the following conditions holds:
\begin{enumroman}
\item A map that is an $I$-cofibration and a weak equivalence is a
$J$-cofibration, or

\item A map that is both a $J$-injective and a weak equivalence is
an $I$-injective.
\end{enumroman}
\end{enumerate}
Then there is a cofibrantly generated model category structure on
$\mathcal M$ in which $I$ is a set of generating cofibrations and
$J$ is a set of generating acyclic cofibrations.
\end{theorem}

We now define our notion of ``equivalence" between two model
categories. Recall that for categories $\mathcal C$ and $\mathcal
D$ a pair of functors
\[ \xymatrix@1{F: \mathcal C \ar@<.5ex>[r] & \mathcal D:R
\ar@<.5ex>[l]} \] is an \emph{adjoint pair} if for each object $X$
of $\mathcal C$ and object $Y$ of $\mathcal D$ there is an
isomorphism $\varphi :\Hom_\mathcal D(FX,Y) \rightarrow
\Hom_\mathcal C(X,RY)$ which is natural in $X$ and $Y$
\cite[IV.1]{macl}.

\begin{definition} \cite[1.3.1]{hovey}
If $\mathcal C$ and $\mathcal D$ are model categories, then the
adjoint pair
\[ \xymatrix@1{F: \mathcal C \ar@<.5ex>[r] & \mathcal D:R
\ar@<.5ex>[l]} \] is a \emph{Quillen pair} if one of the following
equivalent statements is true:
\begin{enumerate}
\item $F$ preserves cofibrations and acyclic cofibrations.

\item $R$ preserves fibrations and acyclic fibrations.
\end{enumerate}
\end{definition}

\begin{definition} \cite[1.3.12]{hovey}
A Quillen pair is a \emph{Quillen equivalence} if for all
cofibrant $X$ in $\mathcal C$ and fibrant $Y$ in $\mathcal D$, a
map $f \colon FX \rightarrow Y$ is a weak equivalence in $\mathcal
D$ if and only if the map $\varphi f \colon X \rightarrow RY$ is a
weak equivalence in $\mathcal C$.
\end{definition}

We will use the following proposition to prove that a Quillen pair
is a Quillen equivalence.  Recall that a functor $F \colon
\mathcal C \rightarrow \mathcal D$ \emph{reflects} a property if,
for any morphism $f$ of $\mathcal C$, whenever $Ff$ has the
property, then so does $f$.

\begin{prop}\cite[1.3.16]{hovey} \label{QPair}
Suppose that
\[ \xymatrix@1{F: \mathcal C \ar@<.5ex>[r] & \mathcal D:R
\ar@<.5ex>[l]} \] is a Quillen pair.  Then the following
statements are equivalent:
\begin{enumerate}
\item This Quillen pair is a Quillen equivalence.

\item $F$ reflects weak equivalences between cofibrant objects
and, for every fibrant $Y$ in $\mathcal D$, the map $F((RY)^c)
\rightarrow Y$ is a weak equivalence.

\item $R$ reflects weak equivalences between fibrant objects and,
for every cofibrant $X$ in $\mathcal C$, the map $X \rightarrow
R((FX)^f)$ is a weak equivalence.
\end{enumerate}
\end{prop}

The existence of a Quillen equivalence between two model
categories is actually a stronger condition than we need, but it
is a convenient way to show that two homotopy theories are the
same. Here, we take the viewpoint that simplicial categories are
models for homotopy theories.  A \emph{simplicial category} is a
category $\mathcal C$ enriched over simplicial sets, or a category
such that, for objects $x$ and $y$ of $\mathcal C$, there is a
simplicial set of morphisms $\Hom_\mathcal C(x,y)$ between them.
We will use the following notion of equivalence of simplicial
categories.

\begin{definition} \cite[2.4]{dk1} \label{DKdefn}
A functor $f \colon \mathcal C \rightarrow \mathcal D$ between two
simplicial categories is a \emph{DK-equivalence} if it satisfies
the following two conditions:
\begin{enumerate}
\item for any objects $x$ and $y$ of $\mathcal C$, the induced map
$\Hom_\mathcal C(x,y) \rightarrow \Hom_\mathcal D(fx,fy)$ is a
weak equivalence of simplicial sets, and

\item the induced map of categories of components $\pi_0 f \colon
\pi_0 \mathcal C \rightarrow \pi_0 \mathcal D$ is an equivalence
of categories.
\end{enumerate}
\end{definition}

Recall that the \emph{category of components} $\pi_0 \mathcal C$
of a simplicial category $\mathcal C$ is the category with the
same objects as $\mathcal C$ and such that
\[ \Hom_{\pi_0 \mathcal C}(x,y) = \pi_0 \Hom_\mathcal C(x,y). \]

Now, the following result tells us that model categories which are
Quillen equivalent actually have equivalent homotopy theories.

\begin{prop} \cite[5.4]{dk1} \label{qedk}
Suppose that $\mathcal C$ and $\mathcal D$ are Quillen equivalent
model categories.  Then the simplicial localizations $\mathcal
{LC}$ and $\mathcal {LD}$ are DK-equivalent.
\end{prop}

\setcounter{subsection}{8}

\subsection{Simplicial Objects}
Recall that a \emph{simplicial set} is a functor $\Deltaop
\rightarrow \mathcal Sets$, where the \emph{cosimplicial category}
${\bf \Delta}$ has as objects the finite ordered sets $[n] = \{0,
\ldots ,n\}$ and as morphisms the order-preserving maps, and
$\Deltaop$ is its opposite category. In particular, for $n \geq
0$, we have $\Delta [n]$, the $n$-simplex, $\dot \Delta [n]$, the
boundary of $\Delta [n]$, and, for $n>0$ and $0 \leq k \leq n$,
$V[n,k]$, which is $\dot \Delta [n]$ with the $k$th face removed
\cite[I.1]{gj}. For any simplicial set $X$, we denote by $X_n$ the
image of $[n]$. There are face maps $d_i:X_n \rightarrow X_{n-1}$
for $0 \leq i \leq n$ and degeneracy maps $s_i:X_n \rightarrow
X_{n+1}$ for $0 \leq i \leq n$, satisfying certain compatibility
conditions \cite[I.1]{gj}.  We denote by $|X|$ the topological
space given by geometric realization of the simplicial set $X$
\cite[I.2]{gj}.

There is a model category structure on simplicial sets in which
the weak equivalences are the maps which become weak homotopy
equivalences of topological spaces after geometric realization
\cite[I.11.3]{gj}. We denote this model category structure by
$\SSets$.   Note in particular that it is cofibrantly generated.
The generating cofibrations are the maps $\dot \Delta [m]
\rightarrow \Delta [m]$ for all $m \geq 0$, and the generating
acyclic cofibrations are the maps $V[m,k] \rightarrow \Delta [m]$
for all $m \geq 1$ and $0 \leq k \leq m$ \cite[3.2.1]{hovey}. This
model category structure is Quillen equivalent to the standard
model category structure on topological spaces
\cite[3.6.7]{hovey}.  In light of this fact, we will sometimes
refer to simplicial sets as ``spaces."

More generally, a simplicial object in a category $\mathcal C$ is
a functor $\Deltaop \rightarrow \mathcal C$ \cite[3.1]{hovey}.  In
particular, a \emph{simplicial space} (or \emph{bisimplicial set})
is a functor $\Deltaop \rightarrow \SSets$ \cite[IV.1]{gj}. Given
a simplicial set $X$, we also use $X$ to denote the constant
simplicial space with the simplicial set $X$ in each degree. By
$X^t$ we denote the simplicial space such that $(X^t)_n$ is the
constant simplicial set $X_n$, or the simplicial set which has the
set $X_n$ in each degree.

Notice, however, that our definition of ``simplicial category" in
this paper is inconsistent with this terminology.  There is a more
general notion of simplicial category by which is meant a
simplicial object in the category of small categories. Such a
simplicial category is a functor $\Deltaop \rightarrow \mathcal
Cat$ where $\mathcal Cat$ is the category with objects the small
categories and morphisms the functors between them. Our definition
of simplicial category coincides with this one when the extra
condition is imposed that the face and degeneracy maps be the
identity map on objects \cite[2.1]{dk1}.

We also require the following additional structure on some of our
model category structures.  A \emph{simplicial model category} is
a model category which is also a simplicial category satisfying
two additional axioms \cite[9.1.6]{hirsch}. (Again, the
terminology is potentially confusing because a simplicial model
category is not a simplicial object in the category of model
categories.) The important part of this structure that we use is
the fact that, given objects $X$ and $Y$ of a simplicial model
category, it makes sense to talk about the \emph{function
complex}, or simplicial set $\Map(X,Y)$.

Given a model category $\mathcal M$, or more generally a category
with weak equivalences, a \emph{homotopy function complex}
$\Map^h(X,Y)$ is a simplicial set which is the morphism space
between $X$ and $Y$ in the simplicial localization $\mathcal{LM}$
\cite[\S 4]{dk1}. If $\mathcal M$ is a simplicial model category,
$X$ is cofibrant in $\mathcal M$, and $Y$ is fibrant in $\mathcal
M$, then $\Map^h(X,Y)$ is weakly equivalent to $\Map(X,Y)$.

\subsection{Localized Model Category Structures}
Several of the model category structures that we use are obtained
by localizing a given model category structure with respect to a
map or a set of maps.  Suppose that $S = \{f \colon A \rightarrow
B\}$ is a set of maps with respect to which we would like to
localize a model category (or category with weak equivalences)
$\mathcal M$. We define an $S$-\emph{local} object $W$ to be an
object of $\mathcal M$ such that for any $f \colon A \rightarrow
B$ in $S$, the induced map on homotopy function complexes
\[ f^* \colon \Map^h(B,W) \rightarrow \Map^h(A,W) \]
is a weak equivalence of simplicial sets.  (If $\mathcal M$ is a
model category, a local object is usually required to be fibrant.)
A map $g \colon X \rightarrow Y$ in $\mathcal M$ is then defined
to be an $S$-\emph{local equivalence} if for every local object
$W$, the induced map on homotopy function complexes
\[ g^* \colon \Map^h(Y,W) \rightarrow \Map^h(X,W) \]
is a weak equivalence of simplicial sets.

The following theorem holds for model categories $\mathcal M$
which are left proper and cellular. We will not define these
conditions here, but refer the reader to \cite[13.1.1,
12.1.1]{hirsch} for more details.  We do note, in particular, that
a cellular model category is cofibrantly generated.  All the model
categories that we localize in this paper can be shown to satisfy
both these conditions.

\setcounter{equation}{10}

\begin{theorem}\cite[4.1.1]{hirsch} \label{Loc}
Let $\mathcal M$ be a left proper cellular model category.  There
is a model category structure $\mathcal L_S \mathcal M$ on the
underlying category of $\mathcal M$ such that:
\begin{enumerate}
\item The weak equivalences are the $S$-local equivalences.

\item The cofibrations are precisely the cofibrations of $\mathcal
M$.

\item The fibrations are the maps which have the right lifting
property with respect to the maps which are both cofibrations and
$S$-local equivalences.

\item The fibrant objects are the $S$-local objects which are
fibrant in $\mathcal M$.

\item If $\mathcal M$ is a simplicial model category, then its
simplicial structure induces a simplicial structure on $\mathcal
L_S \mathcal M$.
\end{enumerate}
\end{theorem}

In particular, given an object $X$ of $\mathcal M$, we can talk
about its functorial fibrant replacement $LX$ in $\mathcal L_S
\mathcal M$. The object $LX$ is an $S$-local object which is
fibrant in $\mathcal M$, and we will refer to it as the
\emph{localization} of $X$ in $\mathcal L_S \mathcal M$.

\setcounter{subsection}{12}

\subsection{Model Category Structures for Diagrams of Spaces}

Suppose that $\mathcal D$ is a small category and consider the
category of functors $\mathcal D \rightarrow \SSets$, denoted
$\SSetsd$.  This category is also called the category of $\mathcal
D$-diagrams of spaces.  We would like to consider model category
structures on $\SSetsd$.

A natural choice for the weak equivalences in $\SSets^\mathcal D$
is the class of levelwise weak equivalences of simplicial sets.
Namely, given two $\mathcal D$-diagrams $X$ and $Y$, we define a
map $f \colon X \rightarrow Y$ to be a weak equivalence if and
only if for each object $d$ of $\mathcal D$, the map $X(d)
\rightarrow Y(d)$ is a weak equivalence of simplicial sets.

There is a model category structure $\SSetsd_f$ on the category of
$\mathcal D$-diagrams with these weak equivalences and in which
the fibrations are given by levelwise fibrations of simplicial
sets.  The cofibrations in $\SSetsd_f$ are then the maps of
simplicial spaces which have the left lifting property with
respect to the maps which are levelwise acyclic fibrations. This
model structure is often called the \emph{projective} model
category structure on $\mathcal D$-diagrams of spaces \cite[IX,
1.4]{gj}. Dually, there is a model category structure $\SSetsd_c$
in which the cofibrations are given by levelwise cofibrations of
simplicial sets, and this model structure is often called the
\emph{injective} model category structure \cite[VIII, 2.4]{gj}.
The small category $\mathcal D$ which we use in this paper is
$\Deltaop$, so that the diagram category $\SSets^{\Deltaop}$ is
just the category of simplicial spaces.

Consider the Reedy model category structure on simplicial spaces
\cite{reedy}. In this structure, the weak equivalences are again
the levelwise weak equivalences of simplicial sets.  The Reedy
model category structure is cofibrantly generated, where the
generating cofibrations are the maps
\[ \dot \Delta[m] \times \Delta [n]^t \cup \Delta [m] \times
\dot \Delta [n]^t \rightarrow \Delta [m] \times \Delta [n]^t \]
for all $n,m \geq 0$.  The generating acyclic cofibrations are the
maps
\[ V[m,k] \times \Delta [n]^t \cup \Delta [m] \times \dot \Delta
[n]^t \rightarrow \Delta [m] \times \Delta [n]^t \] for all $n
\geq 0$, $m \geq 1$, and $0 \leq k \leq m$ \cite[2.4]{rezk}.

It turns out that the Reedy model category structure on simplicial
spaces is exactly the same as the injective model category
structure on this same category, as given by the following result.

\setcounter{equation}{13}

\begin{prop} \cite[15.8.7, 15.8.8]{hirsch} \label{inj}
A map $f \colon X \rightarrow Y$ of simplicial spaces is a
cofibration in the Reedy model category structure if and only if
it is a monomorphism. In particular, every simplicial space is
Reedy cofibrant.
\end{prop}

In light of this result, we denote the Reedy model structure on
simplicial spaces by $\SSets^{\Deltaop}_c$.  Both
$\SSets^{\Deltaop}_c$ and $\SSets^{\Deltaop}_f$ are simplicial
model categories. In each case, given two simplicial spaces $X$
and $Y$, we can define $\Map(X,Y)$ by
\[ \Map (X,Y)_n = \Hom (X \times \Delta [n],Y) \]
where the set on the right-hand side consists of maps of
simplicial spaces.

To establish some notation we need later in the paper, we recall
the definition of fibration in the Reedy model category structure.
If $X$ is a simplicial space, let $\sk_nX$ denote its
$n$-skeleton, generated by the spaces in degrees less than or
equal to $n$, and let $\cosk_nX$ denote the $n$-coskeleton of $X$
\cite[\S 1]{reedy}. A map $X \rightarrow Y$ is a fibration in
$\SSets^{\Deltaop}_c$ if
\begin{itemize}
\item $X_0 \rightarrow Y_0$ is a fibration of simplicial sets, and

\item for all $n \geq 1$, the map $X_n \rightarrow P_n$ is a
fibration, where $P_n$ is defined to be the pullback in the
following diagram:
\[ \xymatrix{P_n \ar[r] \ar[d] & Y_n \ar[d] \\
(\cosk_{n-1}X)_n \ar[r] & (\cosk_{n-1}Y)_n} \]
\end{itemize}

Notice in particular that this pullback diagram is actually a
homotopy pullback diagram, as follows.  If $f \colon X \rightarrow
Y$ is a Reedy fibration, then it has the right lifting property
with respect to all Reedy acyclic cofibrations.  In particular,
there is a dotted arrow lift in the following diagram, where $m
\geq 1$, $0 \leq k \leq m$, and $n \geq 0$:
\[ \xymatrix{V[m,k] \times \dot \Delta [n]^t \ar[r] \ar[d] & X
\ar[d] \\
\Delta [m] \times \dot \Delta [n]^t \ar[r] \ar@{-->}[ur] & Y.} \]

Since the functors $\sk_n$ and $\cosk_n$ are adjoint \cite[\S
1]{reedy}, we have that
\[ (\cosk_{n-1}X)_n \simeq \Map(\Delta [n], \cosk_nX) \simeq \Map(\sk_n \Delta [n], X)
\simeq \Map(\dot \Delta [n], X). \] Therefore, we have a dotted
arrow lift in each diagram
\[ \xymatrix{V[m,k] \ar[r] \ar[d] & (\cosk_{n-1}X)_n \ar[d] \\
\Delta [m] \ar[r] \ar@{-->}[ur] & (\cosk_{n-1}Y)_n.} \]  In
particular, the right-hand vertical arrow is a fibration of
simplicial sets.  Thus, the simplicial set $P_n$ is a homotopy
pullback and therefore homotopy invariant.

We also make use of the projective model category structure
$\SSets^{\Deltaop}_f$ on simplicial spaces. This model category is
also cofibrantly generated; the generating cofibrations are the
maps
\[ \dot \Delta [m] \times \Delta [n]^t \rightarrow \Delta [m]
\times \Delta [n]^t \] for all $m,n \geq 0$ \cite[IV.3.1]{gj}.

In the next section, we localize the Reedy (or injective) and
projective model category structures on simplicial spaces with
respect to a map to obtain model category structures in which the
fibrant objects are Segal spaces (Definition \ref{SegalSpace}). We
will further localize them to obtain model category structures in
which the fibrant objects are complete Segal spaces (Definition
\ref{CSegalSpace}).

\section{Some Definitions and Model Category Structures}

In this section, we define and discuss in turn the three main
structures that we will use in the course of this paper:
simplicial categories, complete Segal spaces, and Segal
categories.

\subsection{Simplicial Categories}
Simplicial categories, most simply stated, are categories enriched
over simplicial sets, or categories with a simplicial set of
morphisms between any two objects.  So, given any objects $x$ and
$y$ in a simplicial category $\mathcal C$, there is a simplicial
set $\Hom_\mathcal C(x,y)$.

Fix an object set $\mathcal O$ and consider the category of
simplicial categories with object set $\mathcal O$ such that all
morphisms are the identity on the objects.  Dwyer and Kan define a
model category structure $\sco$ in which the weak equivalences are
the functors $f \colon \mathcal C \rightarrow \mathcal D$ of
simplicial categories such that given any objects $x$ and $y$ of
$\mathcal C$, the induced map
\[ \Hom_\mathcal C(x,y) \rightarrow \Hom_\mathcal D(x,y) \]
is a weak equivalence of simplicial sets \cite[\S 7]{dk}.  The
fibrations are the functors $f \colon \mathcal C \rightarrow
\mathcal D$ for which these same induced maps are fibrations, and
the cofibrations are the functors which have the left lifting
property with respect to the acyclic fibrations.

It is more useful, however, to consider the category of all small
simplicial categories with no restriction on the objects.  Before
describing the model category structure on this category, we need
a few definitions. Recall from Definition \ref{DKdefn} above that
if $\mathcal C$ is a simplicial category, then we denote by $\pi_0
\mathcal C$ the category of components of $\mathcal C$.

If $\mathcal C$ is a simplicial category and $x$ and $y$ are
objects of $\mathcal C$, a morphism $e \in \Hom_\mathcal C(x,y)_0$
is a \emph{homotopy equivalence} if the image of $e$ in $\pi_0
\mathcal C$ is an isomorphism.

\setcounter{equation}{1}

\begin{theorem} \cite[1.1]{simpcat} \label{SC}
There is a model category structure on the category $\mathcal{SC}$
of small simplicial categories defined by the following three
classes of morphisms:

\begin{enumerate}
\item The weak equivalences are the maps $f \colon \mathcal C
\rightarrow \mathcal D$ satisfying the following two conditions:
\begin{itemize}
\item (W1) For any objects $x$ and $y$ in $\mathcal C$, the map
\[ \Hom_\mathcal C (x,y) \rightarrow \Hom_\mathcal D (fx,fy) \]
is a weak equivalence of simplicial sets.

\item (W2) The induced functor $\pi_0f \colon \pi_0 \mathcal C
\rightarrow \pi_0 \mathcal D$ on the categories of components is
an equivalence of categories.
\end{itemize}

\item The fibrations are the maps $f \colon \mathcal C \rightarrow
\mathcal D$ satisfying the following two conditions:
\begin{itemize}
\item (F1) For any objects $x$ and $y$ in $\mathcal C$, the map
\[ \Hom_\mathcal C (x,y) \rightarrow \Hom_\mathcal D (fx,fy) \]
is a fibration of simplicial sets.

\item (F2) For any object $x_1$ in $\mathcal C$, $y$ in $\mathcal
D$, and homotopy equivalence $e \colon fx_1 \rightarrow y$ in
$\mathcal D$, there is an object $x_2$ in $\mathcal C$ and
homotopy equivalence $d \colon x_1 \rightarrow x_2$ in $\mathcal
C$ such that $fd=e$.
\end{itemize}

\item The cofibrations are the maps which have the left lifting
property with respect to the maps which are fibrations and weak
equivalences.
\end{enumerate}
\end{theorem}

Notice that the weak equivalences are precisely the
DK-equivalences that we defined above (Definition \ref{DKdefn}).

The proof of this theorem actually shows that this model category
structure is cofibrantly generated. Define the functor $U: \SSets
\rightarrow \mathcal{SC}$ such that for any simplicial set $K$,
the simplicial category $UK$ has two objects, $x$ and $y$, and
only nonidentity morphisms the simplicial set $K = \Hom(x,y)$.
Using this functor, we define the generating cofibrations to be
the maps of simplicial categories
\begin{itemize}
\item (C1) $U \dot \Delta [n] \rightarrow U \Delta [n]$ for $n
\geq 0$, and

\item (C2) $\phi \rightarrow \{x\}$, where $\phi$ is the
simplicial category with no objects and $\{x\}$ denotes the
simplicial category with one object $x$ and no nonidentity
morphisms.
\end{itemize}
The generating acyclic cofibrations are defined similarly \cite[\S
1]{simpcat}.

\setcounter{subsection}{2}

\subsection{Segal Spaces and Complete Segal Spaces} \label{sscss}

Complete Segal spaces, defined by Rezk in \cite{rezk}, are more
difficult to describe, but ultimately they are actually easier to
work with than simplicial categories. The name ``Segal" refers to
the similarity between Segal spaces and Segal's $\Gamma$-spaces
\cite{segal}.

We begin by defining Segal spaces.  In \cite[4.1]{rezk}, Rezk
defines for each $0 \leq i \leq k-1$ a map $\alpha^i \colon [1]
\rightarrow [k]$ in ${\bf \Delta}$ such that $0 \mapsto i$ and $1
\mapsto i+1$. Then for each $k$ he defines the simplicial space
\[ G(k)^t= \bigcup_{i=0}^{k-1} \alpha^i \Delta [1]^t \subset \Delta
[k]^t. \]

He shows that, for any simplicial space $X$, there is a weak
equivalence of simplicial sets
\[ \Map^h_{\SSets^{\Deltaop}} (G(k)^t, X) \rightarrow \underbrace{X_1 \times^h_{X_0} \cdots
\times^h_{X_0} X_1}_k, \] where the right hand side is the
homotopy limit of the diagram
\[ \xymatrix{X_1 \ar[r]^{d_0} & X_0 & X_1 \ar[l]_{d_1}
\ar[r]^{d_0} & \ldots \ar[r]^{d_0} & X_0 & X_1 \ar[l]_{d_1}} \]
with $k$ copies of $X_1$.

Now, given any $k$, define the map $\varphi^k \colon G(k)^t
\rightarrow \Delta [k]^t$ to be the inclusion map.  Then for any
simplicial space $W$ there is a map
\[ \varphi_k = \Map^h_{\SSets^{\Deltaop}}(\varphi^k, W) \colon \Map^h_{\SSets^{\Deltaop}}(\Delta
[k]^t,W) \rightarrow \Map^h_{\SSets^{\Deltaop}}(G(k)^t,W). \]
More simply written, this map is
\[ \varphi_k \colon W_k \rightarrow \underbrace{W_1 \times^h_{W_0} \cdots
\times^h_{W_0} W_1}_k \] and is often called a \emph{Segal map}.

\setcounter{equation}{3}

\begin{definition} \cite[4.1]{rezk} \label{SegalSpace}
A Reedy fibrant simplicial space $W$ is a \emph{Segal space} if
for each $k \geq 2$ the map $\varphi_k$ is a weak equivalence of
simplicial sets. In other words, the Segal maps
\[ \varphi_k: W_k \rightarrow \underbrace{W_1 \times^h_{W_0} \cdots
\times^h_{W_0} W_1}_k \] are weak equivalences for all $k \geq 2$.
\end{definition}

Notice that if $W$ is a Segal space, or more generally if $W$ is
Reedy fibrant, we can use ordinary function complexes and a limit
in the definition of the Segal maps \cite[\S 4]{rezk}.

Rezk defines the coproduct of all these inclusion maps
\[ \varphi = \coprod_{k \geq 0} (\varphi^k \colon G(k)^t \rightarrow \Delta [k]^t).
\] Using this map $\varphi$, we have the following result.

\begin{theorem}\cite[7.1]{rezk}
There is a model category structure on simplicial spaces which can
be obtained by localizing the Reedy model category structure with
respect to the map $\varphi$. This model category structure has
the following properties :

\begin{enumerate}
\item The weak equivalences are the maps $f$ for which
$\Map^h_{\SSets^{\Deltaop}}(f,W)$ is a weak equivalence of
simplicial sets for any Segal space $W$.

\item The cofibrations are the monomorphisms.

\item The fibrant objects are the Reedy fibrant $\varphi$-local
objects, which are precisely the Segal spaces.
\end{enumerate}
\end{theorem}

We will refer to this model category structure on simplicial
spaces as the \emph{Segal space model category structure} and
denote it $\sesp_c$.

The properties of Segal spaces enable us to speak of them much in
the same way that we speak of categories.  Heuristically, a simple
example of a Segal space is the nerve of a category $\mathcal C$,
regarded as a simplicial space $\nerve(\mathcal C)^t$. (We need to
take a Reedy fibrant replacement of this nerve to be an actual
Segal space.) In particular, we can define ``objects" and ``maps"
of a Segal space. We summarize the particular details here that we
need; a full description is given by Rezk \cite[\S 5]{rezk}.

Given a Segal space $W$, define its set of \emph{objects}, denoted
$\ob(W)$, to be the set of 0-simplices of the space $W_0$, namely,
the set $W_{0,0}$. Given any two objects $x,y$ in $\ob (W)$,
define the \emph{mapping space} $\map_W(x,y)$ to be the homotopy
fiber of the map $(d_1, d_0) \colon W_1 \rightarrow W_0 \times
W_0$ over $(x,y)$. (Note that since $W$ is Reedy fibrant, this map
is a fibration, and therefore in this case we can just take the
fiber.) Given a 0-simplex $x$ of $W_0$, we denote by $\text{id}_x$
the image of the degeneracy map $s_0 \colon W_0 \rightarrow W_1$.
We say that two 0-simplices of $\map_W(x,y)$, say $f$ and $g$, are
\emph{homotopic}, denoted $f \sim g$, if they lie in the same
component of the simplicial set $\map_W(x,y)$.

Given $f \in \map_W (x,y)_0$ and $g \in \map_W(y,z)_0$, there is a
composite $g \circ f \in \map_W(x,z)_0$, and this notion of
composition is associative up to homotopy.  We define the
\emph{homotopy category} $\Ho(W)$ of $W$ to have as objects the
set $\ob (W)$ and as morphisms between any two objects $x$ and
$y$, the set $\map_{\Ho(W)}(x,y) = \pi_0 \map_W(x,y)$.

A map $g$ in $\map_W(x,y)_0$ is a \emph{homotopy equivalence} if
there exist maps $f,h \in \map_W(y,x)_0$ such that $g \circ f \sim
\text{id}_y$ and $h \circ g \sim \text{id}_x$. Any map in the same
component as a homotopy equivalence is itself a homotopy
equivalence \cite[5.8]{rezk}.  Therefore we can define the space
$W_{\hoequiv}$ to be the subspace of $W_1$ given by the components
whose zero-simplices are homotopy equivalences.

We then note that the degeneracy map $s_0 \colon W_0 \rightarrow
W_1$ factors through $W_\hoequiv$ since for any object $x$ the map
$s_0(x)= \text{id}_x$ is a homotopy equivalence.  Therefore, we
have the following definition:

\begin{definition} \cite[\S 6]{rezk} \label{CSegalSpace}
A \emph{complete Segal space} is a Segal space $W$ for which the
map $s_0 \colon W_0 \rightarrow W_\hoequiv$ is a weak equivalence
of simplicial sets.
\end{definition}

We now consider an alternate way of defining a complete Segal
space which is less intuitive but will enable us to localize the
Segal space model category structure further in such a way that
the complete Segal spaces are the new fibrant objects. Consider
the category $I[1]$ which consists of two objects $x$ and $y$ and
exactly two non-identity maps which are inverse to one another, $x
\rightarrow y$ and $y \rightarrow x$. Denote by $E$ the nerve of
this category, and by $E^t$ the corresponding simplicial space.
There are two maps $\Delta [0]^t \rightarrow E^t$ given by the
inclusions of $\Delta [0]^t$ to the objects $x$ and $y$,
respectively.  Let $\psi \colon \Delta [0]^t \rightarrow E^t$ be
the map which takes $\Delta [0]^t$ to the object $x$.  (It does
not actually matter which one of the two maps we have chosen, as
long as it is fixed.) This map then induces, for any Segal space
$W$, a map on homotopy function complexes
\[ \psi^* \colon \Map^h_{\SSets^{\Deltaop}}(E^t,W) \rightarrow \Map^h_{\SSets^{\Deltaop}}(\Delta [0]^t, W) = W_0. \]

\begin{prop}\cite[6.4]{rezk}
For any Segal space $W$, the map $\psi^*$ of homotopy function
complexes is a weak equivalence of simplicial sets if and only if
$W$ is a complete Segal space.
\end{prop}

Given this proposition, we can further localize the category of
simplicial spaces with respect to this map.

\begin{theorem}\cite[7.2]{rezk} \label{CSS}
Taking the localization of the Reedy model category structure on
simplicial spaces with respect to the maps $\varphi$ and $\psi$
above results in a model category structure which satisfies the
following properties:

\begin{enumerate}
\item The weak equivalences are the maps $f$ such that
$\Map^h_{\SSets^{\Deltaop}}(f,W)$ is a weak equivalence of
simplicial sets for any complete Segal space $W$.

\item The cofibrations are the monomorphisms.

\item The fibrant objects are the complete Segal spaces.
\end{enumerate}
\end{theorem}

We refer to this model category structure on simplicial spaces as
the \emph{complete Segal space model category structure}, denoted
$\css$.  It turns out that when the objects involved are Segal
spaces, the weak equivalences in this model category structure can
be described more explicitly.

\begin{definition}
A map $f \colon U \rightarrow V$ of Segal spaces is a
\emph{DK-equivalence} if
\begin{enumerate}
\item for any pair of objects $x,y \in U_0$, the induced map
$\map_U(x,y) \rightarrow \map_V(fx,fy)$ is a weak equivalence of
simplicial sets, and

\item the induced map $\Ho (f) \colon \Ho(U) \rightarrow \Ho(V)$
is an equivalence of categories.
\end{enumerate}
\end{definition}

We then have the following result by Rezk:

\begin{theorem}\cite[7.7]{rezk} \label{SSDK}
Let $f \colon U \rightarrow V$ be a map of Segal spaces.  Then $f$
is a DK-equivalence if and only if it becomes a weak equivalence
in $\mathcal{CSS}$.
\end{theorem}

Note that these weak equivalences have been given the same name as
the ones in $\mathcal{SC}$.  While this may at first seem strange,
the two definitions are very similar, in fact rely on the same
generalization of the idea of equivalence of categories to a
simplicial setting.

However, what is especially nice about the complete Segal space
model category structure is the simple characterization of the
weak equivalences between the fibrant objects.

\begin{prop} \cite[7.6]{rezk} \label{DKReedy}
A map $f \colon U \rightarrow V$ between complete Segal spaces is
a DK-equivalence if and only if it is a levelwise weak
equivalence.
\end{prop}

This proposition is actually a special case of a more general
result.  In any localized model category structure, a map is a
local equivalence between fibrant objects if and only if it is a
weak equivalence in the original model category structure
\cite[3.2.18]{hirsch}.

It is also possible to localize the projective model category
structure $\SSets^{\Deltaop}_f$ on the category of simplicial
spaces to obtain analogous model category structures. We will
denote the localization of the projective model category structure
by with respect to the map $\varphi$ by $\sesp_f$. There is also a
localization of the projective model category structure with
respect to the maps $\varphi$ and $\psi$ analogous to the model
category structure $\mathcal{CSS}$, but we do not need this
structure here.

\setcounter{subsection}{11}

\subsection{Segal Categories} \label{secats}
Lastly, we consider the Segal categories.  We begin by defining
the preliminary notion of a Segal precategory.

\setcounter{equation}{12}

\begin{definition} \cite[\S 2]{hs}
A \emph{Segal precategory} is a simplicial space $X$ such that the
simplicial set $X_0$ in degree zero is discrete, i.e., a constant
simplicial set.
\end{definition}

In the case of Segal precategories, it again makes sense to talk
about the Segal maps
\[\varphi_k \colon X_k \rightarrow \underbrace{X_1 \times^h_{X_0} \cdots \times^h_{X_0}
X_1}_k \] for each $k \geq 2$. Since $X_0$ is discrete, we can
actually take the limit
\[ \underbrace{X_1 \times_{X_0} \cdots \times_{X_0} X_1}_k \] on the
right-hand side.

\begin{definition} \cite[\S 2]{hs} \label{SegalCat}
A \emph{Segal category} $X$ is a Segal precategory such that each
Segal map $\varphi_k$ is a weak equivalence of simplicial sets for
$k \geq 2$.
\end{definition}

Note that the definition of a Segal category is similar to that of
a Segal space, with the additional requirement that the degree
zero space be discrete.  However, Segal categories are not
required to be Reedy fibrant, so they are not necessarily Segal
spaces.

Given a fixed set $\mathcal O$, we can consider the category
$\SSets^{\Deltaop}_\mathcal O$ whose objects are the Segal
precategories with $\mathcal O$ in degree zero and whose morphisms
are the identity on this set. There is a model category structure
$\SSets^{\Deltaop}_{\mathcal O, f}$ on this category in which the
weak equivalences are levelwise \cite[3.7]{simpmon}. In other
words, $f \colon X \rightarrow Y$ is a weak equivalence if for
each $n \geq 0$, the map $f_n \colon X_n \rightarrow Y_n$ is a
weak equivalence of simplicial sets. Furthermore, the fibrations
are also levelwise. This model structure can then be localized
with respect to a map similar to the map which we used to obtain
the Segal space model category structure.

We first need to determine what this map should be.  We begin by
considering the maps of simplicial spaces $\varphi^k \colon G(k)^t
\rightarrow \Delta [k]^t$ and adapting them to the case at hand.

The first problem is that $\Delta [k]^t$ is not going to be in
$\SSets^{\Deltaop}_{\mathcal O, f}$ for all values of $k$.
Instead, we need to define a separate $k$-simplex for any
$k$-tuple $x_0, \ldots ,x_k$ of objects in $\mathcal O$, denoted
$\Delta [k]^t_{x_0, \ldots ,x_k}$, so that the objects are
preserved. Note that this object $\Delta [k]^t_{x_0, \ldots ,x_k}$
also needs to have all elements of $\mathcal O$ as 0-simplices, so
we add any of these elements that have not already been included
in the $x_i$'s, plus their degeneracies in higher degrees.

Then we can define
\[ G(k)^t_{x_0, \ldots ,x_k} = \bigcup_{i=0}^{k-1} \alpha^i \Delta
[1]^t_{x_i, x_{i+1}}. \] Now, we need to take coproducts not only
over all values of $k$, but also over all $k$-tuples of vertices.
Hence, the resulting map $\varphi_\mathcal O$ looks like
\[ \varphi_\mathcal O = \coprod_{k \geq 0}(\coprod_{(x_0, \ldots ,x_k) \in
\mathcal O^{k+1}} (G(k)^t_{x_0, \ldots ,x_k} \rightarrow \Delta
[k]^t_{x_0, \ldots ,x_k})). \] Setting $\xu =(x_0, \ldots, x_k)$,
we can write the component maps as $G(k)^t_\xu \rightarrow \Delta
[k]^t_\xu$.  We can then localize $\SSets^{\Deltaop}_{\mathcal
O,f}$ with respect to the map $\varphi_\mathcal O$ to obtain a
model category which we denote $\mathcal
{LSS}ets^{\Deltaop}_{\mathcal O,f}$.

There are also analogous model category structures
$\SSets^{\Deltaop}_{\mathcal O, c}$ and
$\mathcal{LSS}ets^{\Deltaop}_{\mathcal O,c}$ on the category of
Segal precategories with a fixed set $\mathcal O$ in degree zero
with the same weak equivalences but where the cofibrations, rather
than the fibrations, are defined levelwise, and then we can
localize with respect to the same map \cite[3.9]{simpmon},
\cite[A.1.1]{tv}.

However, we would like a model category structure on the category
of all Segal precategories, not just on these more restrictive
subcategories.  In the course of this paper, we prove the
existence of two model category structures on Segal precategories.
Unlike in the fixed object set case, we cannot actually obtain the
model category structure via localization of a model category
structure with levelwise weak equivalences since it is not
possible to put a model structure on the category of Segal
precategories in which the weak equivalences are levelwise and in
which the cofibrations are monomorphisms.

To see that there is no such model structure, suppose that one did
exist and consider the map $f \colon \Delta [0]^t \amalg \Delta
[0]^t \rightarrow \Delta [0]^t$. By model category axiom MC5, $f$
could be factored as the composite of a cofibration $\Delta [0]^t
\amalg \Delta [0]^t \rightarrow X$ followed by an acyclic
fibration $X \rightarrow \Delta [0]^t$.  However, since the weak
equivalences would be levelwise weak equivalences, $X_0$ would
have to consist of one point. However, the only map $(\Delta [0]^t
\amalg \Delta [0]^t)_0 \rightarrow X_0$ is not a monomorphism.
Thus, there is no such factorization of the map $f$, and therefore
there can be no model category structure satisfying the two given
properties.

\setcounter{subsection}{14}

\subsection{Relationship Between Simplicial Categories and Segal
Categories in Fixed Object Set Cases} \label{scosecato}

Recall from above that there is a model category structure $\sco$
on the category whose objects are the simplicial categories with a
fixed set $\mathcal O$ of objects and whose morphisms are the
functors which are the identity on the objects and that there is a
model category structure $\mathcal{LSS}ets^{\Deltaop}_{\mathcal O,
f}$ on the category whose objects are the Segal precategories with
the set $\mathcal O$ in degree zero and whose morphisms are the
identity on degree zero.

\setcounter{equation}{15}

\begin{theorem} \cite[5.5]{simpmon} \label{Oequiv}
There is an adjoint pair
\[ \xymatrix@1{F_\mathcal O: \mathcal{LSS}ets^{\Deltaop}_{\mathcal O, f} \ar@<.5ex>[r] & \sco:
R_\mathcal O \ar@<.5ex>[l]} \] which is a Quillen equivalence.
\end{theorem}

The proof of this theorem uses a generalization of a result by
Badzioch \cite[6.5]{bad} which relates strict and homotopy
algebras over an algebraic theory.  This generalization uses the
notion of multi-sorted algebraic theory \cite{multisort}.

A key step in this proof is a explicit description of the
localization of the objects $\Delta [n]^t_\xu$. Up to homotopy,
this localization is the same as the localization of the objects
$G(n)^t_\xu$ and is obtained by taking the colimit of stages of a
filtration
\[ G(n)^t_\xu = \Psi_1 G(n)^t_\xu \subseteq \Psi_2 G(n)^t_\xu
\subseteq \cdots \]

Let $e_i$ denote the nondegenerate 1-simplex $x_{i-1} \rightarrow
x_i$ in $G(n)^t_\xu$ and let $w_j$ denote a word in the $e_i$'s
which can be obtained via ``composition" of these 1-simplices. The
$k$-th stage of the filtration is given by
\[ (\Psi_k (G(n)^t_\xu))_m=\left\{(w_1 \mid \cdots \mid w_m) \mid \ell(w_1 \cdots w_m) \leq k \right\} \]
where $\ell(w_1 \cdots w_n)$ denotes the length of the word $w_1
\cdots w_n$. The colimit of this filtration is weakly equivalent
to $L_cG(n)^t_\xu$ in $\LSSets^{\Deltaop}_{\mathcal O,f}$.

We show in the proof of \cite[4.2]{simpmon} that for each $i \geq
1$ the map
\[ \Psi_i G(n)^t_\xu \rightarrow \Psi_{i+1} G(n)^t_\xu \] is a
DK-equivalence, and that the unique map from $G(n)^t_\xu$ to the
colimit of this directed system is also a DK-equivalence.

In the current paper, we use some of the ideas of the proof from
the fixed object set case, but we no longer use multi-sorted
theories as we pass from $\sco$ to $\mathcal {SC}$ and
$\SSets^{\Deltaop}_\mathcal O$ to $\Secat$.

\section{Methods of Obtaining Segal Precategories from Simplicial
Spaces}

In the course of proving the existence of these two model category
structures $\Secat_c$ and $\Secat_f$, we need sets of generating
cofibrations and generating acyclic cofibrations which are similar
to those of the Reedy and projective model category structures on
simplicial spaces. However, we need to modify these maps so that
they are actually maps between Segal precategories.  The purpose
of this section is to define two methods of modifying the
generating cofibrations and generating acyclic cofibrations so
that they are actually maps between Segal precategories, and to
prove a result which we need to prove the existence of the model
structures $\Secat_c$ and $\Secat_f$.

The first method we call reduction, and we use it to define the
generating cofibrations in $\Secat_c$. Consider the forgetful
functor from the category of Segal precategories to the category
of simplicial spaces. This map has a left adjoint, which we call
the reduction map. Given a simplicial space $X$, we denote its
reduction by $(X)_r$.  The degree $n$ space of $(X)_r$ is obtained
from $X_n$ by collapsing the subspace $s_0^nX_0$ of $X_n$ to the
discrete space $\pi_0 (s_0^n X_0)$, where $s_0^n$ is the iterated
degeneracy map.

Recall that the cofibrations in the Reedy model category structure
on simplicial spaces are monomorphisms (Proposition \ref{inj}) and
that the Reedy generating cofibrations are of the form
\[ \dot \Delta [m] \times \Delta [n]^t \cup \Delta [m] \times \dot
\Delta [n]^t \rightarrow \Delta [m] \times \Delta [n]^t \] for all
$n,m \geq 0$.  In general, these maps are not in $\Secat$ because
the objects involved are not Segal precategories. Therefore, we
apply this reduction functor to these maps.

Thus, we consider the maps
\[ (\dot \Delta [m] \times \Delta [n]^t \cup \Delta [m] \times
\dot \Delta [n]^t)_r \rightarrow (\Delta [m] \times \Delta
[n]^t)_r. \] However, we still need to make some modifications to
assure that all these maps are actually monomorphisms.  In
particular, we need to check the case where $n=0$.  If $n=m=0$,
and if $\phi$ denotes the empty simplicial space, we obtain the
map $\phi \rightarrow \Delta [0]^t$, which is a monomorphism.
However, when $n=0$ and $m=1$, we get the map $\Delta [0]^t \amalg
\Delta [0]^t \rightarrow \Delta [0]^t$, which is not a
monomorphism. When $n=0$ and $m \geq 2$, we obtain the map $\Delta
[0]^t \rightarrow \Delta [0]^t$.  This map is an isomorphism, and
thus there is no reason to include it in the generating set.
Therefore, we define the set
\[ I_c= \{(\dot \Delta
[m] \times \Delta [n]^t \cup \Delta [m] \times \dot \Delta
[n]^t)_r \rightarrow (\Delta [m] \times \Delta [n]^t)_r \}
\] for all $m \geq 0$ when $n \geq 1$ and for $n=m=0$.  This set $I_c$ will be a set of
generating cofibrations of $\Secat_c$.

This reduction process works well in almost all situations, but we
have problems when we try to reduce some of the generating
cofibrations in $\SSets^{\Deltaop}_f$, namely the maps
\[ \dot \Delta [1] \times \Delta [n]^t \rightarrow \Delta [1]
\times \Delta [n]^t \] for any $n \geq 0$.  The object $\Delta [1]
\times \Delta [n]^t$ reduces to a Segal precategory with $n+1$
points in degree zero, but the object $\dot \Delta [1] \times
\Delta [n]^t$ reduces to a Segal precategory with $2(n+1)$ points
in degree zero.  In other words, the reduced map in this case is
no longer a monomorphism.

Consider the set $\Delta [n]_0$ and denote by $\Delta [n]^t_0$ the
doubly constant simplicial space defined by it.  For $m \geq 1$
and $n \geq 0$, define $P_{m,n}$ to be the pushout of the diagram
\[ \xymatrix{\dot \Delta [m] \times \Delta[n]^t_0 \ar[r] \ar[d] &
\dot \Delta [m] \times \Delta [n]^t \ar[d] \\
\Delta [n]^t_0 \ar[r] & P_{m,n}.} \]  If $m=0$, then we define
$P_{m,0}$ to be the empty simplicial space.  For all $m \geq 0$
and $n \geq 1$, define $Q_{m,n}$ to be the pushout of the diagram
\[ \xymatrix{\Delta [m] \times \Delta [n]^t_0 \ar[r] \ar[d] &
\Delta [m] \times \Delta [n]^t \ar[d] \\
\Delta [n]^t_0 \ar[r] & Q_{m,n}.} \] For each $m$ and $n$, the map
$\dot \Delta [m] \times \Delta [n]^t$ induces a map $i_{m,n}
\colon P_{m,n} \rightarrow Q_{m,n}$.  We then define the set $I_f
= \{i_{m,n} \colon P_{m,n} \rightarrow Q_{m,n} \mid m,n \geq 0\}$.
Note that when $m \geq 2$ this construction gives exactly the same
objects as those given by reduction, namely that $P_{m,n}$ is
precisely $(\dot \Delta [m] \times \Delta [n]^t)_r$ and likewise
$Q_{m,n}$ is precisely $(\Delta [m] \times \Delta [n]^t)_r$.

Given a Segal precategory $X$, we denote by $X_n(v_0, \ldots v_n)$
the fiber of the map $X_n \rightarrow X_0^{n+1}$ over $(v_0,
\ldots ,v_n) \in X_0^{n+1}$, where this map is given by iterated
face maps of $X$.  More specifically, $X_0^{n+1}= (\cosk_0X)_n$
and the map $X_n \rightarrow X_0^{n+1}$ is given by the map $X
\rightarrow \cosk_0X$.

If $\Hom$ denotes morphism set and $X$ is an arbitrary simplicial
space, notice that we can use the pushout diagrams defining the
objects $P_{m,n}$ and $Q_{m,n}$ to see that
\[ \Hom(P_{m,n},X) \cong \coprod_{v_0, \ldots ,v_n} \Hom(\dot \Delta [m],
X_n(v_0, \ldots ,v_n)) \] and
\[ \Hom(Q_{m,n},X) \cong \coprod_{v_0, \ldots ,v_n} \Hom(\Delta [m],
X_n(v_0, \ldots ,v_n)). \]

We now state and prove a lemma using the maps in $I_f$.

\begin{lemma} \label{PQ}
Suppose a map $f \colon X \rightarrow Y$ has the right lifting
property with respect to the maps in $I_f$.  Then the map $X_0
\rightarrow Y_0$ is surjective and each map
\[ X_n(v_0, \ldots ,v_n) \rightarrow Y_n(fv_0, \ldots ,fv_n) \] is
an acyclic fibration of simplicial sets for each $n \geq 1$ and
$(v_0, \ldots ,v_n) \in X_0^{n+1}$.
\end{lemma}

\begin{proof}
The surjectivity of $X_0 \rightarrow Y_0$ follows from the fact
that $f$ has the right lifting property with respect to the map
$P_{0,0} \rightarrow Q_{0,0}$.

In order to prove the remaining statement, it suffices to show
that there is a dotted arrow lift in any diagram of the form
\begin{equation} \label{mn}
\xymatrix{\dot \Delta [m] \ar[r] \ar[d] & X_n(v_0, \ldots ,v_n)
\ar[d] \\
\Delta [m] \ar[r] \ar@{-->}[ur] & Y_n(fv_0, \ldots ,fv_n)}
\end{equation}
for $m,n \geq 0$.

By our hypothesis, there is a dotted arrow lift in diagrams of the
form
\begin{equation} \label{pq}
\xymatrix{P_{m,n} \ar[r] \ar[d] & X \ar[d] \\
Q_{m,n} \ar[r] \ar@{-->}[ur] & Y}
\end{equation}
for all $m,n \geq 0$. The existence of the lift in diagram
\ref{pq} is equivalent to the surjectivity of the map
$\Hom(Q_{m,n},X) \rightarrow P$ in the following diagram, where
$P$ denotes the pullback and $\Hom$ denotes morphism set:
\[ \xymatrix{\Hom(Q_{m,n},X) \ar[r] & P \ar[r] \ar[d] & \Hom(P_{m,n},X) \ar[d] \\
& \Hom(Q_{m,n},Y) \ar[r] & \Hom(P_{m,n},Y).} \]

Now, as noted above we have that
\[ \Hom(Q_{m,n},X) \cong \coprod_{v_0, \ldots ,v_n} \Hom(\Delta [m],
X_n(v_0, \ldots ,v_n)) \] and analogous weak equivalences for the
other objects of the diagram.

Using these weak equivalences and being particularly careful in
the cases where $m=1$ and $m=0$, one can show that for each $m,n
\geq 0$ the dotted-arrow lift in diagram \ref{mn} exists and
therefore that each map
\[ X_n(v_0, \ldots ,v_n) \rightarrow Y_n(fv_0, \ldots ,fv_n) \] is
an acyclic fibration of simplicial sets for each $n \geq 1$.
\end{proof}

\section{A Segal Category Model Category Structure on Segal Precategories}

In this section, we prove the existence of the model category
structure $\Secat_c$.

We would like to define a functorial ``localization" functor $L_c$
on $\Secat$ such that, given any Segal precategory $X$, its
localization $L_cX$ is a Segal space which is a Segal category
weakly equivalent to $X$ in $\sesp_c$. We begin by considering a
functorial localization functor in $\sesp_c$ and then modifying it
so that it takes values in $\Secat$.  In the case of $\sesp_c$,
this localization functor is precisely the functorial fibrant
replacement functor.

A choice of generating acyclic cofibrations for $\sesp_c$ is the
set of maps
\[ V[m,k] \times \Delta [n]^t \cup \Delta [m] \times G(n)^t \rightarrow \Delta [m] \times \Delta [n]^t \]
for $n \geq 0$, $m \geq 1$, and $0 \leq k \leq m$ \cite[\S
4.2]{hirsch}.  Therefore, one can use the small object argument to
construct a functorial localization functor by taking a colimit of
pushouts, each of which is along the coproduct of all these maps
\cite[\S 4.3]{hirsch}.

If we apply this functor to a Segal precategory, the maps with
$n=0$ are problematic because taking pushouts along them will not
result in a space which is discrete in degree zero. We claim that
we can obtain a functorial localization functor $L_c$ on the
category $\Secat$ by taking a colimit of iterated pushouts along
the maps
\[ V[m,k] \times \Delta [n]^t \cup \Delta [m] \times G(n)^t \rightarrow \Delta [m] \times \Delta [n]^t
\]
for $n,m \geq 1$ and $0 \leq k \leq m$.

To see that this restricted set of maps is sufficient, consider a
Segal precategory $X$ and the Segal category $L_cX$ we obtain from
taking such a colimit. Then for any $0 \leq k \leq m$, consider
the diagram
\[ \xymatrix{V[m,k] \ar[r] \ar[d] & \Map^h(G(0)^t, L_cX) \ar[d] \\
\Delta [m] \ar[r] \ar@{-->}[ur] & \Map^h(\Delta [0]^t, L_cX).} \]
Since $\Delta [0]^t$ is isomorphic to $G(0)^t$, and since $L_cX$
is discrete in degree zero, the right-hand vertical map is an
isomorphism of discrete simplicial sets.  Therefore, a dotted
arrow lift exists in this diagram.  It follows that the map $L_cX
\rightarrow \Delta [0]^t$ has the right lifting property with
respect to the maps
\[ V[m,k] \times \Delta [n]^t \cup \Delta [m] \times G(n)^t \rightarrow \Delta [m] \times \Delta [n]^t \]
for all $n \geq 0$, $m \geq 1$, and $0 \leq k \leq m$.  Therefore,
$L_cX$ is fibrant in $\sesp_c$, namely, a Segal space.

Since $L_cX$ is a Segal space, it makes sense to talk about the
mapping space $\map_{L_cX}(x,y)$ and the homotopy category $\Ho
(L_cX)$.  Given these facts, we show that there exists a model
category structure $\Secat_c$ on Segal precategories with the
following three distinguished classes of morphisms:

\begin{itemize}
\item Weak equivalences are the maps $f:X \rightarrow Y$ such that
the induced map $L_cX \rightarrow L_cY$ is a DK-equivalence of
Segal spaces. (Again, we will call such maps
\emph{DK-equivalences}.)

\item Cofibrations are the monomorphisms.  (In particular, every
Segal precategory is cofibrant.)

\item Fibrations are the maps with the right lifting property with
respect to the maps which are both cofibrations and weak
equivalences.
\end{itemize}

\begin{theorem} \label{Secatc}
There is a cofibrantly generated model category structure
$\Secat_c$ on the category of Segal precategories with the above
weak equivalences, fibrations, and cofibrations.
\end{theorem}

We first need to define sets $I_c$ and $J_c$ as our candidates for
generating cofibrations and generating acyclic cofibrations,
respectively.

We take as generating cofibrations the set \[ I_c= \{(\dot \Delta
[m] \times \Delta [n]^t \cup \Delta [m] \times \dot \Delta
[n]^t)_r \rightarrow (\Delta [m] \times \Delta [n]^t)_r \}
\] for all $m \geq 0$ when $n \geq 1$ and for $n=m=0$.  Notice
that since taking a pushout along such a map amounts to attaching
an $m$-simplex to the space in degree $n$, any cofibration can be
written as a directed colimit of pushouts along the maps of $I_c$.

We then define the set $J_c=\{i \colon A \rightarrow B\}$ to be a
set of representatives of isomorphism classes of maps in $\Secat$
satisfying two conditions:
\begin{enumerate}
\item For all $n \geq 0$, the spaces $A_n$ and $B_n$ have
countably many simplices.

\item The map $i \colon A \rightarrow B$ is a monomorphism and a
weak equivalence.
\end{enumerate}

Given these proposed generating acyclic cofibrations, we need to
show that any acyclic cofibration in $\Secat_c$ is a directed
colimit of pushouts along such maps.  To prove this result, we
require several lemmas. The proofs of the first three we omit
here; proofs can be found in the author's thesis \cite{thesis}.

\begin{lemma} \label{cone}
Let $A \rightarrow B$ be a CW-inclusion.  The following statements
are equivalent:
\begin{enumerate}
\item $A \rightarrow B$ is a weak equivalence of topological
spaces.

\item For all $n \geq 1$, any map of pairs $(D^n,S^{n-1})
\rightarrow (B,A)$ extends over the map of cones
$(CD^n,CS^{n-1})$.

\item For all $n \geq 1$, any map $(D^n, S^{n-1}) \rightarrow
(B,A)$ is homotopic to a constant map.
\end{enumerate}
\end{lemma}

\begin{lemma} \label{CWpair}
Let $f \colon X \rightarrow Y$ be a an inclusion of simplicial
sets which is a weak equivalence, and let $W$ and $Z$ be
simplicial sets such that we have a diagram of inclusions
\[ \xymatrix{W \ar[r] \ar[d] & Z \ar[d] \\
X \ar[r] & Y} \] Let $u \colon (D^n,S^{n-1}) \rightarrow
(|Z|,|W|)$ be a relative map of CW-pairs.  Then the inclusion $i
\colon (|Z|,|W|) \rightarrow (|Y|,|X|)$ can be factored as a
composite
\[ (|Z|,|W|) \rightarrow (|K|,|L|) \rightarrow (|Y|,|X|) \] where
$K$ is a subspace of $Y$ obtained from $Z$ by attaching a finite
number of nondegenerate simplices, $L$ is a subspace of $X$, and
the composite map of relative CW-complexes
\[ (D^n,S^{n-1}) \rightarrow (|Z|,|W|) \rightarrow (|K|,|L|) \] is
homotopic rel $S^{n-1}$ to a map $D^n \rightarrow |L|$.
\end{lemma}

\begin{lemma} \label{countable}
Let $(Y,X)$ be a CW-pair such that $X$ and $Y$ have only countably
many cells.  Then for a fixed $n \geq 0$, there are only countably
many homotopy classes of maps $(D^n,S^{n-1}) \rightarrow (Y,X)$.
\end{lemma}

If $A \rightarrow B$ is a monomorphism of Segal precategories,
then taking the localization via the small object argument gives
us that $L_cA \rightarrow L_cB$ is a monomorphism of Segal
categories.  In particular, if $A \subseteq B$ is an inclusion,
then we can regard $L_cA \subseteq L_cB$ as an inclusion as well.

\begin{lemma} \label{finite}
Let $A$ and $B$ be Segal precategories such that $A \subseteq B$.
Let $\sigma$ be a simplex in $L_cB$ which is not in $L_cA$.  Then
there exists a Segal precategory $A'$ such that $A'$ is obtained
from $A$ by attaching a finite number of nondegenerate simplices
and $\sigma$ is in $LA'$.
\end{lemma}

\begin{proof}
By our description of our localization functor at the beginning of
the section, $L_cB$ is obtained from $B$ by taking a colimit of
pushouts, each of which is along the map
\[ \coprod_{m, k, n} V[m,k] \times \Delta [n]^t \cup \Delta [m] \times G(n)^t
\rightarrow \coprod_{m, k, n} \Delta [m] \times \Delta [n]^t \]
for $n,m \geq 1$ and $0 \leq k \leq m$.  The Segal category $L_cB$
is the colimit of a filtration
\[ B \subseteq \Psi^1 B \subseteq \Psi^2 B \subseteq \cdots \] where
each $\Psi^i$ is given by a colimit of iterated pushouts along
this map. Since $\sigma$ is a single simplex, it is small and
therefore $\sigma$ is in $\Psi^n B$ for some $n$.

Therefore, $\sigma$ is obtained by attaching $\Delta [m] \times
\Delta [n]^t$ along a finite number of nondegenerate simplices of
$\Psi^{n-1}B$.  We can then apply the preceding argument to each
of these simplices and inductively obtain a finite number of
nondegenerate simplices of $B$ which form a sub-Segal precategory
which we will call $C$.  We then define $A'=A \cup C$.
\end{proof}

We then state one more lemma, which is a generalization of a lemma
given by Hirschhorn \cite[2.3.6]{hirsch}.

\begin{lemma} \label{Commute}
Let the map $g \colon A \rightarrow B$ be an inclusion of Segal
precategories, each of which has countably many simplices. If $X$
is a Segal precategory with countably many simplices, then its
localization $LX$ with respect to the map $g$ has only countably
many simplices.
\end{lemma}

We are now able to state and prove our result about generating
cofibrations.

\begin{prop} \label{Cofs}
Any acyclic cofibration $j \colon C \rightarrow D$ in $\Secat_c$
can be written as a directed colimit of pushouts along the maps in
$J_c$.
\end{prop}

\begin{proof}
Note that by definition $j \colon C \rightarrow D$ is a
monomorphism of Segal precategories.  We assume that it is an
inclusion.  Let $U$ be a subsimplicial space of $D$ such that $U$
has countably many simplices in each degree. Apply the
localization functor $L_c$ to obtain a diagram of Segal categories
\[ \xymatrix{L_c(U \cap C) \ar[r] \ar[d] & L_cU \ar[d] \\
L_cC \ar[r] & L_cD.} \]  Since $U$ has only countably many
simplices, this localization process adds at most a countable
number of simplices to the original simplicial space by Lemma
\ref{Commute}.

We would like to find a Segal precategory $W$ such that $U
\subseteq W \subseteq D$ and such that the map $W \cap C
\rightarrow W$ is in the set $J_c$.

First consider the map
\[ \Ho(L_c(U \cap C)) \rightarrow \Ho(L_cU) \] which we want to be an
equivalence of categories.  If it is not an equivalence, then
there exists $z \in (L_cU)_0$ which is not equivalent to some $z'
\in (L_c(U \cap C))_0$. However, there is such a $z'$ when we
consider $z$ as an element of $(L_cD)_0$, since $j:C \rightarrow
D$ is a DK-equivalence. If this $z'$ is not in $(U \cap C)_0$,
then we add it.  Repeat this process for all such $z$.

Now for each such $z$, consider the four mapping spaces in $L_cU$
involving the objects $z$ and $z'$: $\map_{L_cU}(z,z)$,
$\map_{L_cU}(z,z')$, $\map_{L_cU}(z',z)$, and
$\map_{L_cU}(z',z')$. We want the sets of components of these four
spaces to be isomorphic to one another in $\Ho(L_cU)$. We can
attach a countable number of simplices via an analogous argument
to the one in the proof of Lemma \ref{finite} such that these sets
of components are isomorphic. We then repeat the same argument to
assure that $\pi_0\map_{L_cU}(x,z)$ is isomorphic to
$\pi_0\map_{L_cU}(x,z')$ for each $x \in U_0$ and analogously for
the sets of components of the mapping spaces out of each such $x$.

By repeating this process for each such $z$, we obtain a Segal
precategory $Y$ with a countable number of simplices such that
$\Ho(L_c(Y \cap C)) \rightarrow \Ho(L_cY)$ is an equivalence of
categories. However, we do not necessarily have that for each $x,y
\in L_c(Y \cap C)$,
\[ \map_{L_c(Y \cap C)}(x,y) \rightarrow \map_{L_cY}(x,y) \] is a
weak equivalence of simplicial sets.  Therefore we consider all
maps
\[ (D^n,S^{n-1}) \rightarrow
(|\map_{L_cY}(x,y)|,|\map_{L_c(Y \cap C)}(x,y)|) \rightarrow
(|\map_{L_cD}(x,y)|,|\map_{L_cC}(x,y)|) \] for each $x,y \in (Y
\cap C)_0$ and $n \geq 0$.  Identify all $x,y$, and $n$ such that
the map
\[ (D^n,S^{n-1}) \rightarrow
(|\map_{L_cY}(x,y)|,|\map_{L_c(Y \cap C)}(x,y)|) \] is not
homotopic to a constant map.

However each composite map
\[ (D^n,S^{n-1}) \rightarrow
(|\map_{L_cY}(x,y)|,|\map_{L_c(Y \cap C)}(x,y)|) \rightarrow
(|\map_{L_cD}(x,y)|,|\map_{L_cC}(x,y)|) \] is homotopic to a
constant map by Lemma \ref{cone} since
\[ |\map_{L_cC}(x,y)| \rightarrow |\map_{L_cD}(x,y)| \] is a weak
equivalence.

For each such $x$, $y$, and $n$, it follows from Lemma
\ref{CWpair} that there exists some pair of simplicial sets
\[(\map_{L_cY}(x,y), \map_{L_c(Y \cap C)}(x,y)) \subseteq (K,L) \subseteq
(\map_{L_cD}(x,y), \map_{L_cC}(x,y)) \] such that the composite
map
\[ (D^n,S^{n-1}) \rightarrow (|\map_{L_cY}(x,y)|,|\map_{L_c(Y \cap C)}(x,y)|) \rightarrow
(|K|,|L|) \] is homotopic to a constant map, and the pair $(K,L)$
is obtained from the pair $(\map_{L_cY}(x,y),\map_{L_c(Y \cap
C)}(x,y))$ by attaching a finite number of nondegenerate
simplices.  We apply Lemma \ref{finite} to each of these new
simplices obtained by considering each nontrivial homotopy class
to obtain some Segal precategory $Y'$ with a countable number of
number of simplices such that each composite map
\[ (D^n,S^{n-1}) \rightarrow (|\map_{L_cY}(x,y)|,|\map_{L_c(Y \cap C)}(x,y)|) \rightarrow
(|\map_{L_cY'}(x,y)|,|\map_{L_c(Y' \cap C)}(x,y)|) \] is homotopic
to a constant map.

However, the process of adding simplices may have created more
maps
\[ (D^n,S^{n-1}) \rightarrow (|\map_{L_cY'}(x,y)|, |\map_{L_c(Y'
\cap C)}(x,y)|) \] that are not homotopic to a constant map.
Therefore we repeat this argument, perhaps countably many times,
until, taking a colimit over all of them, we obtain a Segal
precategory $W$ such that each map
\[ (D^n,S^{n-1}) \rightarrow (|\map_{L_cW}(x,y)|,|\map_{L_c(W \cap
C)}(x,y)|)\] is homotopic to a constant map. Since each of these
steps added only countably many simplices to the original Segal
precategory $U$, and since by Lemma \ref{cone}
\[ \map_{L_c(W \cap C)}(x,y) \rightarrow \map_{L_cW}(x,y) \] is a
weak equivalence for all $x,y \in (L_c(W \cap C))_0$, the map $W
\cap C \rightarrow W$ is in the set $J_c$.

Now, take some $\widetilde U$ obtained from $W$ by adding a
countable number of simplices, consider the inclusion map
$\widetilde U \cap C \rightarrow \widetilde U$, and repeat the
entire process. To show that we can repeat this argument, taking a
(possibly transfinite) colimit, and eventually obtain the map $j:C
\rightarrow D$, it suffices to show that the localization functor
$L_c$ commutes with arbitrary directed colimits of inclusions.
However, this fact follows from \cite[2.2.18]{hirsch}.
\end{proof}

Now, we have two definitions of acyclic fibration that we need to
show coincide: the fibrations which are weak equivalences, and the
maps with the right lifting property with respect to the maps in
$I_c$.

\begin{prop} \label{RLPC}
The maps with the right lifting property with respect to the maps
in $I_c$ are fibrations and weak equivalences.
\end{prop}

Before giving a proof of this proposition, we begin by looking at
the maps in $I_c$ and determining what an $I_c$-injective looks
like. Recall the definition of the coskeleton of a simplicial
space from the paragraph following Proposition \ref{inj}.  If $f
\colon X \rightarrow Y$ has the right lifting property with
respect to the maps in $I_c$, then for each $n \geq 1$, the map
$X_n \rightarrow P_n$ is an acyclic fibration of simplicial sets,
where $P_n$ is the pullback in the diagram
\[ \xymatrix{P_n \ar[r] \ar[d] & Y_n \ar[d] \\
(\cosk_{n-1}X)_n \ar[r] & (\cosk_{n-1}Y)_n.} \] In the case that
$n=0$, the restrictions on $m$ and $n$ give us that the map $X_0
\rightarrow Y_0$ is a surjection rather than the isomorphism we
get in the Reedy case. Notice that by the same argument given for
the Reedy model category structure (in the section following
Proposition \ref{inj} above), the simplicial sets $P_n$ can be
characterized up to weak equivalence as homotopy pullbacks and are
therefore homotopy invariant.

This characterization of the maps with the right lifting property
with respect to $I_c$ will enable us to prove Proposition
\ref{RLPC}. Before proceeding to the proof, however, we state a
lemma, whose proof we defer to section 9.

\begin{lemma} \label{xnyn}
Suppose that $f:X \rightarrow Y$ is a map of Segal precategories
which is an $I_c$-injective. Then $f$ is a DK-equivalence.
\end{lemma}

\begin{proof} [Proof of Proposition \ref{RLPC}]
Suppose that $f \colon X \rightarrow Y$ is an $I_c$-injective, or
a map which has the right lifting property with respect to the
maps in $I_c$. Note that $f$ then has the right lifting property
with respect to all cofibrations. Since, in particular, it has the
right lifting property with respect to the acyclic cofibrations,
it is a fibration by definition.  It remains to show that $f$ is a
weak equivalence.

However, this fact follows from Lemma \ref{xnyn}, proving the
proposition.
\end{proof}

We now state the converse, which we prove in section 9.

\begin{prop} \label{FibWE}
The maps in $\Secat_c$ which are both fibrations and weak
equivalences are $I_c$-injectives.
\end{prop}

Now we prove a lemma which we need to check the last condition for
our model category structure.

\begin{lemma} \label{PushoutJ}
A pushout along a map of $J_c$ is also an acyclic cofibration in
$\Secat_c$.
\end{lemma}

\begin{proof}
Let $j \colon A \rightarrow B$ be a map in $J_c$.  Notice that $j$
is an acyclic cofibration in the model category $\css$.  Since
$\css$ is a model category, we know that a pushout along an
acyclic cofibration is again an acyclic cofibration
\cite[3.14(ii)]{ds}. If all the objects involved are Segal
precategories, then the pushout will again be a Segal precategory
and therefore the pushout map will be an acyclic cofibration in
$\Secat_c$.
\end{proof}

\begin{prop} \label{Jcofc}
If a map of Segal precategories is a $J_c$-cofibration, then it is
an $I_c$-cofibration and a weak equivalence.
\end{prop}

\begin{proof}
By definition and Proposition \ref{Cofs}, a $J_c$-cofibration is a
map with the left lifting property with respect to the maps with
the right lifting property with respect to the acyclic
cofibrations.  However, by the definition of fibration, these maps
are the ones with the left lifting property with respect to the
fibrations.

Similarly, using Propositions \ref{RLPC} and \ref{FibWE}, an
$I_c$-cofibration is a map with the left lifting property with
respect to the acyclic fibrations.  Thus, we need to show that a
map with the left lifting property with respect to the fibrations
has the left lifting property with respect to the acyclic
fibrations and is a weak equivalence.  Since the acyclic
fibrations are fibrations, it remains to show that the maps with
the left lifting property with respect to the fibrations are weak
equivalences.

Let $f \colon A \rightarrow B$ be a map with the left lifting
property with respect to all fibrations.  By Lemma \ref{PushoutJ}
above, we know that a pushout along maps of $J_c$ is an acyclic
cofibration. Therefore, we can use the small object argument
\cite[10.5.15]{hirsch} to factor the map $f \colon A \rightarrow
B$ as the composite of an acyclic cofibration $A \rightarrow A'$
and a fibration $A' \rightarrow B$. Then there exists a dotted
arrow lift in the diagram
\[ \xymatrix{A \ar[r]^\simeq \ar[d] & A' \ar[d] \\
B \ar[r]^{\text{id}} \ar@{-->}[ur] & B} \] showing that the map $A
\rightarrow B$ is a retract of the map $A \rightarrow A'$ and
therefore a weak equivalence.
\end{proof}

\begin{proof}[Proof of Theorem \ref{Secatc}]
Axiom MC1 follows since limits and colimits of Segal precategories
(computed as simplicial spaces) still have discrete zero space and
are therefore Segal precategories.  MC2 and MC3 (for weak
equivalences) work as usual, for example see \cite[8.10]{ds}.

It remains to show that the four conditions of Theorem
\ref{CofGen} are satisfied.  The set $I_c$ permits the small
object argument because the generating cofibrations in the Reedy
model category structure do.  We can show that the objects $A$
which appear as the sources of the maps in $J_c$ are small using
an analogous argument to the one for simplicial sets
\cite[10.4.4]{hirsch}, so the set $J_c$ permits the small object
argument. Thus, condition 1 is satisfied.

Condition 2 is precisely the statement of Proposition \ref{Jcofc}.
Condition 3 and condition 4(ii) are precisely the statements of
Propositions \ref{RLPC} and \ref{FibWE}.
\end{proof}

Note that the reduced Reedy acyclic cofibrations
\[ (V[m,k] \times \Delta [n]^t \cup \Delta [m] \times \dot \Delta [n]^t)_r
\rightarrow (\Delta [m] \times \Delta [n]^t)_r \] are acyclic
cofibrations in $\Secat_c$ for $m \geq 0$ when $n \geq 1$ and for
$n=m=0$.

\begin{cor}
The fibrant objects in $\Secat_c$ are Reedy fibrant Segal
categories.
\end{cor}

\begin{proof}
Suppose that $X$ is fibrant in $\Secat_c$.  Then, since the
reduced Reedy cofibrations are acyclic cofibrations in $\Secat_c$
and since $X$ has discrete zero space, it follows that $X$ is
Reedy fibrant.

Then, since the maps
\[ (\Delta [m] \times G(n)^t)_r \rightarrow (\Delta [m] \times
\Delta [n]^t)_r \] for all $m,n \geq 0$ are acyclic cofibrations
in $\Secat_c$, it follows that $X$ is a Segal category.
\end{proof}

The converse statement, that the Reedy fibrant Segal categories
are fibrant in $\Secat_c$, also holds \cite{fibrant}.

\section{A Quillen Equivalence Between $\Secat_c$ and $\mathcal{CSS}$}

In this section, we will show that there is a Quillen equivalence
between the model category structure $\Secat_c$ on Segal
precategories and the complete Segal space model category
structure $\css$ on simplicial spaces.  We first need to show that
we have an adjoint pair of maps between the two categories.

Let $I \colon \Secat_c \rightarrow \css$ be the inclusion functor
of Segal precategories into the category of all simplicial spaces.
We will show that there is a right adjoint functor $R \colon \css
\rightarrow \Secat_c$ which ``discretizes" the degree zero space.

Let $W$ be a simplicial space.  Define simplicial spaces $U=
\cosk_0(W_0)$ and $V= \cosk_0(W_{0,0})$. There exist maps $W
\rightarrow U \leftarrow V$. Then we take the pullback $RW$ in the
diagram
\[ \xymatrix{RW \ar[r] \ar[d] & V \ar[d] \\
W \ar[r] & U.} \]

Note that $RW$ is a Segal precategory.  If $W$ is a complete Segal
space, then so are $U$ and $V$, and in this case $RW$ is a Segal
category, which we can see as follows.  The pullback at degree 1
gives
\[ \xymatrix{(RW)_1 \ar[r] \ar[d] & W_{0,0} \times W_{0,0} \ar[d]
\\
W_1 \ar[r] & W_0 \times W_0} \] and at degree 2 we get
\[ \xymatrix{(RW)_1 \times_{(RW)_0} (RW)_1 \ar[r] \ar[d] &
(W_{0,0})^3 \ar[d] \\
W_2 \simeq W_1 \times_{W_0} W_1 \ar[r] & W_0 \times W_0 \times
W_0.} \]  Looking at these pullbacks, and the analogous ones for
higher $n$, we notice that $RW$ is in fact a Segal category.

We define the functor $R \colon \css \rightarrow \Secat_c$ which
takes a simplicial space $W$ to the Segal precategory $RW$ given
by the description above.

\begin{prop}
The functor $R \colon \css \rightarrow \Secat_c$ is right adjoint
to the inclusion map $I \colon \Secat_c \rightarrow \css$.
\end{prop}

\begin{proof}
We need to show that there is an isomorphism
\[ \Hom_{\Secat_c} (Y,RW) \cong  \Hom_{\css}(IY, W) \]
for any Segal precategory $Y$ and simplicial space $W$.

Suppose that we have a map $Y = IY \rightarrow W$.  Since $Y$ is a
Segal precategory, $Y_0$ is equal to $Y_{0,0}$ viewed as a
constant simplicial set.  Therefore, we can restrict this map to a
unique map $Y \rightarrow V$, where $V$ is the Segal precategory
defined above. Then, given the universal property of pullbacks,
there is a unique map $Y \rightarrow RW$. Hence, we obtain a map
\[ \varphi \colon \Hom_{\css}(IY,W) \rightarrow \Hom_{\Secat_c}(Y,RW). \]

This map is surjective because given any map $Y \rightarrow RW$ we
can compose it with the map $RW \rightarrow W$ to obtain a map $Y
\rightarrow W$.

Now for any Segal precategory $Y$, consider the diagram
\[ \xymatrix{Y \ar@/^/[rrd] \ar@/_/[ddr] \ar@{-->}[dr] && \\
& RW \ar[r] \ar[d] & V \ar[d] \\
& W \ar[r] & U} \] Because this diagram must commute and the image
of the map $Y_0 \rightarrow W_0$ is contained in $W_{0,0}$ since
$Y$ is a Segal precategory, this map uniquely determines what the
map $Y \rightarrow V$ has to be.  Therefore, given a map $Y
\rightarrow RW$, it could only have come from one map $Y
\rightarrow W$. Thus, $\varphi$ is injective.
\end{proof}

Now, we need to show that this adjoint pair respects the model
category structures that we have.

\begin{prop}
The adjoint pair of functors
\[ \xymatrix@1{I:\Secat_c \ar@<.5ex>[r] & \css :R \ar@<.5ex>[l]} \] is
a Quillen pair.
\end{prop}

\begin{proof}
It suffices to show that the inclusion map $I$ preserves
cofibrations and acyclic cofibrations.  $I$ preserves cofibrations
because they are defined to be monomorphisms in each category.
Also in each of the two categories, a map is a weak equivalence if
it is a DK-equivalence after localizing to obtain a Segal space,
as given in Theorem \ref{SSDK}. In each case an acyclic
cofibration is an inclusion satisfying this property.  Therefore,
the map $I$ preserves acyclic cofibrations.
\end{proof}

\begin{theorem}
The Quillen pair
\[ \xymatrix@1{I:\Secat_c \ar@<.5ex>[r] & \css :R \ar@<.5ex>[l]} \] is
a Quillen equivalence.
\end{theorem}

\begin{proof}
We need to show that $I$ reflects weak equivalences between
cofibrant objects and that for any fibrant object $W$ (i.e.,
complete Segal space) in $\css$, the map $I((RW)^c) = IRW
\rightarrow W$ is a weak equivalence in $\Secat_c$.

The fact that $I$ reflects weak equivalences between cofibrant
objects follows from the same argument as the one in the proof of
the Quillen pair.  To prove the second part, it remains to show
that the map $j \colon RW \rightarrow W$ in the pullback diagram
\[ \xymatrix{RW \ar[r] \ar[d]^j & V \ar[d] \\
W \ar[r] & U} \] is a DK-equivalence.  It suffices to show that
the map of objects $\ob(RW) \rightarrow \ob(W)$ is surjective and
that the map $\map_{RW}(x,y) \rightarrow \map_W(jx,jy)$ is a weak
equivalence, where the object set of a Segal space is defined as
in section \ref{sscss}. However, notice by the definition of $RW$
that $\ob(RW)=\ob(W)$. In particular, $jx=x$ and $jy=y$. Then
notice, using the pullback that defines $(RW)_1$, that
$\map_{RW}(x,y) \simeq \map_W(x,y)$. Therefore, the map $RW
\rightarrow W$ is a DK-equivalence.
\end{proof}

\section{Another Segal Category Model Category Structure on Segal Precategories}

The model category structure $\Secat_c$ that we defined above is
helpful for the Quillen equivalence with the complete Segal space
model category structure, but there does not appear to be a
Quillen equivalence between it and the model category structure
$\mathcal{SC}$ on simplicial categories. Therefore, we need
another model category structure $\Secat_f$ to obtain such a
Quillen equivalence.

In the model category structure $\Secat_c$, we started with the
generating cofibrations in the Reedy model category structure and
adapted them to be generating cofibrations of Segal precategories.
In this second model category structure, we use modified
generating cofibrations from the projective model category
structure on simplicial spaces so that the objects involved are
Segal precategories.

We make the following definitions for a model category structure
$\Secat_f$ on the category of Segal precategories.

\begin{itemize}
\item The weak equivalences are the same as those of $\Secat_c$.

\item The cofibrations are the maps which can be formed by taking
iterated pushouts along the maps of the set $I_f$ defined in
section 4.

\item The fibrations are the maps with the right lifting property
with respect to the maps which are cofibrations and weak
equivalences.
\end{itemize}

Notice that to define the weak equivalences in this case we want
to use a functorial localization in $\sesp_f$ rather than
$\sesp_c$. We define a localization functor $L_f$ in the same way
that we defined $L_c$ at the beginning of section 5 but making
necessary changes in light of the fact that we are starting from
the model structure $\sesp_f$.  So, in a sense, the weak
equivalences are not defined identically in the two categories,
since they make use of the same localization of different model
category structures on the category of simplicial spaces. However,
in each case the weak equivalences are the same in the unlocalized
model category, so we can define homotopy function complexes using
only the underlying category and the weak equivalences.  Recall by
the definition of local objects that a map $X \rightarrow Y$ is a
local equivalence if and only if the induced map of homotopy
function complexes
\[ \Map^h(Y,Z) \rightarrow \Map^h(X,Z) \] is a weak equivalence of
simplicial sets for any local $Z$.  In particular, the weak
equivalences of the localized category depend only on the weak
equivalences of the unlocalized category.  Therefore the weak
equivalences in $\Secat_c$ and $\Secat_f$ are actually the same.

\begin{theorem} \label{secatf}
There is a cofibrantly generated model category structure
$\Secat_f$ on the category of Segal precategories in which the
weak equivalences, fibrations, and cofibrations are defined as
above.
\end{theorem}

We define the set $J_f$ to be a set of isomorphism classes of maps
$\{i \colon A \rightarrow B\}$ such that
\begin{enumerate}
\item for all $n \geq 0$, the spaces $A_n$ and $B_n$ have
countably many simplices, and

\item $i:A \rightarrow B$ is an acyclic cofibration.
\end{enumerate}

We would like to show that $I_f$ (defined in section 4) is a set
of generating cofibrations and that $J_f$ is a set of generating
acyclic cofibrations for $\Secat_f$.

We begin with the following lemma.

\begin{lemma} \label{acofs}
Any acyclic cofibration $j \colon C \rightarrow D$ in $\Secat_f$
can be written as a directed colimit of pushouts along the maps in
$J_f$.
\end{lemma}

\begin{proof}
The argument that we used to prove Proposition \ref{Cofs} still
holds, applying the functor $L_f$ rather than $L_c$.
\end{proof}

\begin{prop} \label{AFibs}
A map $f \colon X \rightarrow Y$ is an acyclic fibration in
$\Secat_f$ if and only if it is an $I_f$-injective.
\end{prop}

\begin{proof}
First suppose that $f$ has the right lifting property with respect
to the maps in $I_f$.  Then we claim that for each $n \geq 0$ and
$(v_0, \ldots ,v_n) \in X_0^{n+1}$, the map $X_n(v_0, \ldots ,v_n)
\rightarrow Y_n(fv_0, \ldots ,fv_n)$ is an acyclic fibration of
simplicial sets.  This fact, however, follows from Lemma \ref{PQ}.
In particular, it is a weak equivalence, and therefore we can
apply the proof of Lemma \ref{xnyn} to show that the map $X
\rightarrow Y$ is a DK-equivalence, completing the proof of the
first direction.  (The proof does not follow precisely in this
case, in particular because not all monomorphisms are
cofibrations.  However, we can use the fact that weak equivalences
are the same in $\Secat_c$ and $\Secat_f$ to see that the argument
still holds.)

Then, to prove the converse, assume that $f$ is a fibration and a
weak equivalence.  Then we can apply the proof of Proposition
\ref{FibWE}, making the factorizations in the projective model
category structure rather than in the Reedy model category
structure.  The argument follows analogously.
\end{proof}

\begin{prop} \label{ICofs}
A map in $\Secat_f$ is a $J_f$-cofibration if and only if it is an
$I_f$-cofibration and a weak equivalence.
\end{prop}

\begin{proof}
This proof follows just as the proof of Proposition \ref{Jcofc},
again using the projective structure rather than the Reedy
structure.
\end{proof}

\begin{proof}[Proof of Theorem \ref{secatf}]
As before, we must check the conditions of Theorem \ref{CofGen}.
Condition 1 follows just as in the proof of Theorem \ref{Secatc}.
Condition 2 is precisely the statement of Proposition \ref{ICofs}.
Condition 3 and condition 4(ii) follow from Proposition
\ref{AFibs} after applying Lemma \ref{acofs}.
\end{proof}

We now prove that both our model category structures on the
category of Segal precategories are Quillen equivalent.

\begin{theorem}
The identity functor induces a Quillen equivalence
\[ \xymatrix@1{I:\Secat_f \ar@<.5ex>[r] & \Secat_c:J. \ar@<.5ex>[l]} \]
\end{theorem}

\begin{proof}
Since both maps are the identity functor, they form an adjoint
pair.  We then show that this adjoint pair is a Quillen pair.

We first make some observations between the two categories. Notice
that the cofibrations of $\Secat_f$ form a subclass of the
cofibrations of $\Secat_c$ since they are monomorphisms.
Similarly, the acyclic cofibrations of $\Secat_f$ form a subclass
of the acyclic cofibrations of $\Secat_c$.

In particular, these observations imply that the left adjoint $I
\colon \Secat_f \rightarrow \Secat_c$ preserves cofibrations and
acyclic cofibrations. Hence, we have a Quillen pair.

It remains to show that this Quillen pair is a Quillen
equivalence. To do so, we must show that given any cofibrant $X$
in $\Secat_f$ and fibrant $Y$ in $\Secat_c$, a map $f \colon IX
\rightarrow Y$ is a weak equivalence in $\Secat_f$ if and only if
$\varphi f \colon X \rightarrow JY$ is a weak equivalence in
$\Secat_c$. However, this follows from the fact that the weak
equivalences are the same in each category.
\end{proof}

\begin{note}
One might ask at this point why we could not just use the
$\Secat_f$ model category structure and show a Quillen equivalence
between it and the model category structure $\mathcal{CSS}_f$
where we localize the projective model category structure (rather
than the Reedy) with respect to the maps $\varphi$ and $\psi$. The
existence of such a Quillen equivalence would certainly simplify
this paper!

However, if we work with ``complete Segal spaces" which are
fibrant in the projective model structure rather than in the Reedy
structure, then for a fibrant object $W$ the map $W \rightarrow U$
used in defining the right adjoint $\mathcal{CSS} \rightarrow
\Secat_c$ is no longer necessarily a fibration.  Therefore, the
pullback $RW$ is no longer a homotopy pullback and in particular
not homotopy invariant. If $RW$ is not homotopy invariant, then
there is no guarantee that the map $RW \rightarrow W$ is a
DK-equivalence, and the argument for a Quillen equivalence fails.
Thus, the $\Secat_c$ and $\mathcal{CSS}$ model structures are
necessary.
\end{note}

\section{A Quillen Equivalence Between $\mathcal {SC}$ and
$\Secat_f$}

We begin, as above, by defining an adjoint pair of functors
between the two categories $\mathcal{SC}$ and $\Secat_f$.  We have
the nerve functor $R \colon \mathcal{SC} \rightarrow \Secat_f$.
In order to define a left adjoint to this map, we need some
terminology.

\begin{definition}
Let $\mathcal D$ be a small category and $\SSetsd$ the category of
functors $\mathcal D \rightarrow \SSets$.  Let $S$ be a set of
morphisms in $\SSetsd$.  An object $Y$ of $\SSetsd$ is
\emph{strictly} $S$-\emph{local} if for every morphism $f \colon A
\rightarrow B$ in $S$, the induced map on function complexes
\[ f^* \colon \Map (B,Y) \rightarrow \Map (A,Y) \]
is an isomorphism of simplicial sets. A map $g \colon C
\rightarrow D$ in $\SSetsd$ is a \emph{strict} $S$-\emph{local
equivalence} if for every strictly $S$-local object $Y$ in
$\SSetsd$, the induced map
\[ g^* \colon \Map(D,Y) \rightarrow \Map(C,Y) \]
is an isomorphism of simplicial sets.
\end{definition}

Now, we can view Segal precategories as functors $\Deltaop
\rightarrow \SSets$.  Because we require the image of [0] to be a
discrete simplicial set, the category of Segal precategories is a
subcategory of the category of all such functors.  In this
section, we are going to regard simplicial categories as the
strictly local objects in $\Secat_f$ with respect to the map
$\varphi$ described in section \ref{sscss}.

Although we are actually working in a subcategory, we can still
use the following lemma to obtain a left adjoint functor $F$ to
our inclusion map $R$, since the construction will always produce
a simplicial space with discrete 0-space when applied to such a
simplicial space.

\begin{lemma} \cite[5.6]{multisort}
Consider two categories, the category of all diagrams $X \colon
\mathcal D \rightarrow \SSets$ and the category of strictly local
diagrams with respect to the set of maps $S= \{f \colon A
\rightarrow B\}$. The forgetful functor from the category of
strictly local diagrams to the category of all diagrams has a left
adjoint.
\end{lemma}

We define the functor $F\colon \Secat_f \rightarrow \mathcal {SC}$
to be this left adjoint to the inclusion functor of strictly local
diagrams into all diagrams $R\colon \mathcal {SC} \rightarrow
\Secat_f$.

\begin{prop}
The adjoint pair
\[ \xymatrix@1{F:\Secat_f \ar@<.5ex>[r] & \mathcal{SC} :R \ar@<.5ex>[l]} \] is
a Quillen pair.
\end{prop}

\begin{proof}
We prove that this adjoint pair is a Quillen pair by showing that
the left adjoint $F$ preserves cofibrations and acyclic
cofibrations.  We begin by considering cofibrations.

Since $F$ is a left adjoint functor, it preserves colimits.
Therefore, it suffices to show that $F$ preserves the set $I_f$ of
generating cofibrations in $\Secat_f$.  Recall that the elements
of this set are the maps $P_{m,n} \rightarrow Q_{m,n}$ as defined
in section 4.  We begin by considering the maps $P_{n,1}
\rightarrow Q_{n,1}$ for any $n \geq 0$.  The strict localization
of such a map is precisely the map of simplicial categories $U
\dot \Delta [n] \rightarrow U \Delta [n]$ (section 3.1) which is a
generating cofibration in $\mathcal{SC}$.  We can also see that
the strict localization of any $P_{m,n} \rightarrow Q_{m,n}$ can
be obtained as the colimit of iterated pushouts along the
generating cofibrations of $\mathcal {SC}$.  Therefore, $F$
preserves cofibrations.

We now need to show that $F$ preserves acyclic cofibrations.  To
do so, first consider the model category structure
$\mathcal{LSS}ets^{\Deltaop}_{\mathcal O, f}$ (defined in section
\ref{scosecato}) on Segal precategories with a fixed set $\mathcal
O$ in degree zero and the model category structure $\sco$ of
simplicial categories with a fixed object set $\mathcal O$. Recall
from section \ref{scosecato} that there is a Quillen equivalence
\[ \xymatrix@1{F_\mathcal O:\mathcal{LSS}ets^{\Deltaop}_{\mathcal O, f}
\ar@<.5ex>[r] & \sco :R_\mathcal O. \ar@<.5ex>[l]} \]
In particular, if $X$ is a cofibrant object of
$\mathcal{LSS}ets^{\Deltaop}_{\mathcal O, f}$, then there is a
weak equivalence $X \rightarrow R_\mathcal O((F_\mathcal OX)^f)$.
Notice that $F_\mathcal O$ agrees with $F$ on Segal precategories
with the set $\mathcal O$ in degree zero, and similarly for
$R_\mathcal O$ and $R$.

Suppose, then, that $X$ is an object of
$\mathcal{LSS}ets^{\Deltaop}_{\mathcal O, f}$, $Y$ is an object of
$\mathcal{LSS}ets^{\Deltaop}_{\mathcal O',f}$, and $X \rightarrow
Y$ is an acyclic cofibration in $\Secat_f$.  We have a commutative
diagram
\[ \xymatrix{X \ar[r]^-\simeq \ar[d] & L_fX \ar[d] \\
Y \ar[r]^-\simeq & L_fY} \] where the upper and lower horizontal
maps are weak equivalences not only in $\Secat_f$, but in
$\mathcal{LSS}ets^{\Deltaop}_{\mathcal O, f}$ and
$\mathcal{LSS}ets^{\Deltaop}_{\mathcal O',f}$, respectively.
However, using the fixed-object case Quillen equivalence, the
functors $F_\mathcal O$ and $F_{\mathcal O'}$ (and hence $F$) will
preserve these weak equivalences, giving us a diagram
\[ \xymatrix{FX \ar[r]^-\simeq \ar[d] & FL_fX \ar[d] \\
FY \ar[r]^-\simeq & FL_fX.} \]

Using these weak equivalences and our assumption that $L_fX
\rightarrow L_fY$ is a DK-equivalence, we obtain a diagram
\[ \xymatrix{L_fX \ar[r]^-\simeq \ar[d]^\simeq & RFL_fX \ar[d] \\
L_fY \ar[r]^-\simeq & RFL_fY} \] in which the upper horizontal
arrow is a weak equivalence in $\LSSets^{\Deltaop}_{\mathcal O,f}$
and the lower horizontal arrow is a weak equivalence in
$\LSSets^{\Deltaop}_{\mathcal O',f}$.  The commutativity of this
diagram implies that the map $RFL_fX \rightarrow RFL_fY$ is a
DK-equivalence also. Thus, we have shown that $F$ preserves
acyclic cofibrations between cofibrant objects.

It remains to show that $F$ preserves all acyclic cofibrations.
Suppose that $f \colon X \rightarrow Y$ is an acyclic cofibration
in $\Secat_f$.  Apply the cofibrant replacement functor to the map
$X \rightarrow Y$ to obtain an acyclic cofibration $X' \rightarrow
Y'$, and notice that in the resulting commutative diagram
\[ \xymatrix{X' \ar[r] \ar[d] & Y' \ar[d] \\
X \ar[r] & Y} \] the vertical arrows are levelwise weak
equivalences.

Now consider the following diagram, where the top square is a
pushout diagram:
\[ \xymatrix{X' \ar[r]^\simeq \ar[d] & Y' \ar[d] \\
X \ar[r]^\simeq \ar[d]_= & Y'' \ar[d] \\
X \ar[r]^\simeq & Y.} \]  Notice that all three of the horizontal
arrows are acyclic cofibrations in $\Secat_f$, the upper and lower
by assumption and the middle one because pushouts preserve acyclic
cofibrations \cite[3.14]{ds}. Now we apply the functor $F$ to this
diagram to obtain a diagram
\begin{equation} \label{fxfy}
\xymatrix{FX' \ar[r]^\simeq \ar[d] & FY' \ar[d] \\
FX \ar[d] \ar[r]^\simeq & FY'' \ar[d] \\
FX \ar[r] & FY.}
\end{equation}
The top horizontal arrow is an acyclic cofibration since $F$
preserves acyclic cofibrations between cofibrant objects.
Furthermore, since $F$ is a left adjoint and hence preserves
colimits, the middle horizontal arrow is also an acyclic
cofibration because the top square is a pushout square.

Now, recall that, given an object $X$ in a model category
$\mathcal C$, the \emph{category of objects under} $X$ has as
objects the morphisms $X \rightarrow Y$ in $\mathcal C$ for any
object $Y$, and as morphisms the maps $Y \rightarrow Y'$ in
$\mathcal C$ making the appropriate triangular diagram commute
\cite[7.6.1]{hirsch}.  There is a model category structure on this
under category in which a morphism is a weak equivalence,
fibration, or cofibration if it is in $\mathcal C$
\cite[7.6.5]{hirsch}. In particular, a object $X \rightarrow Y$ is
cofibrant in the under category if it is a cofibration in
$\mathcal C$.

With this definition in mind, to show that the bottom horizontal
arrow of diagram \ref{fxfy} is an acyclic cofibration, consider
the following diagram in the category of cofibrant objects under
$X$:
\[ \xymatrix{X \ar[r] \ar[dr] & Y'' \ar[d] \\
& Y.} \]  Now, let $\mathcal O''$ denote the set in degree zero of
$Y''$ (and also of $Y$) which is not in the image of the map from
$X$.  Now we have the diagram in the category of cofibrant objects
under $X \amalg \mathcal O''$ with the same set in degree zero
\[ \xymatrix{X \amalg \mathcal O'' \ar[r] \ar[dr] & Y'' \ar[d] \\
& Y.} \]  However, since we are now working in a fixed object set
situation, we know by Theorem \ref{Oequiv} that $F_{\mathcal O''}$
is the left adjoint of a Quillen pair, and therefore the map
$F_{\mathcal O''} Y'' \rightarrow F_{\mathcal O''}Y$ is a weak
equivalence in $\mathcal {SC}_{\mathcal O''}$, and in particular a
DK-equivalence when regarded as a map in $\mathcal {SC}$. It
follows that the map $FX \rightarrow FY$ is a weak equivalence,
and $F$ preserves acyclic cofibrations.
\end{proof}

Recall that we are regarding a Segal category as a local diagram
and a simplicial category as a strictly local diagram in
$\Secat_f$.

\begin{lemma} \label{localize}
The map $X \rightarrow FX$ is a DK-equivalence for every cofibrant
object $X$ in $\Secat_f$.
\end{lemma}

\begin{proof}
First consider a free diagram in $\Secat_f$, namely some $\amalg_i
Q_{m_i, n_i}$, where each $Q_{m_i, n_i}$ is defined as in section
4. If $Y$ is a fibrant object in $\Secat_f$, then we have
\[ \begin{aligned}
\Map_{\Secat_f}(\coprod_i Q_{m_i,n_i}, Y) & \simeq \prod_i
\Map_{\Secat_f}(Q_{m_i,n_i},Y) \\
& \simeq \prod_i \coprod_{v_0, \ldots ,v_n} \Map_{\SSets}(\Delta
[m_i], Y_{n_i}(v_0, \ldots,v_n)) \\
& \simeq \prod_i, \coprod_{v_0, \ldots ,v_n} \Map_{\SSets}(\Delta
[0], Y_{n_i}(v_0, \ldots,v_n)) \\
& \simeq \Map_{\Secat_f}(\coprod_i Q_{0,n_i},Y) \\
& \simeq \Map_{\Secat_f}(\coprod_i \Delta [n_i]^t, Y)
\end{aligned} \]
Therefore, it suffices to consider free diagrams $\amalg_i \Delta
[n_i]^t$.  Such a diagram is a Segal category.  It is also the
nerve of a category and thus a strictly local diagram.  It follows
that the map
\[ \coprod_i \Delta [n_i]^t \rightarrow F(\coprod_i \Delta
[n_i]^t) \] is a DK-equivalence.

Now suppose that $X$ is any cofibrant object in $\Secat_f$.  Then
$X$ can be written as a directed colimit $X \simeq
\colim_{\Deltaop} X_j$, where each $X_j$ can be written as
$\amalg_i \Delta[n_i]^t$. As before we regard $FX$ as a strictly
local object in $\Secat_f$.  If $Y$ is a fibrant object in
$\Secat_f$ which is strictly local, we have
\[ \begin{aligned}
\Map_{\Secat_f}(\colim_{\Deltaop}X_j, Y) & \simeq \lim_{\bf
\Delta} \Map_{\Secat_f}(X_j, Y) \\
& \simeq \lim_{\bf \Delta} \Map_{\Secat_f}(FX_j,Y) \\
& \simeq \Map_{\Secat_f}(\colim_{\Deltaop} FX_j,Y) \\
& \simeq \Map_{\Secat_f}(F (\colim_{\Deltaop} (FX_j)),Y)
\end{aligned} \]
We can now apply the result that
\[ F (\colim (FX_j)) \simeq F (\colim X_j). \]  (This fact is proved
in \cite[4.1]{simpmon} for ordinary localization, but it holds for
strict localization in this case since each $X_j$ is cofibrant and
$F$ preserves cofibrant objects.) Therefore we have
\[ \Map_{\Secat_f}(F (\colim_{\Deltaop} (FX_j)),Y) \simeq
\Map_{\Secat_f}(FX,Y). \] It follows that the map $X \rightarrow
FX$ is a DK-equivalence.
\end{proof}

We are now able to prove the main result of this section.

\begin{theorem}
The Quillen pair
\[ \xymatrix@1{F:\Secat_f \ar@<.5ex>[r] & \mathcal{SC} :R \ar@<.5ex>[l]} \] is
a Quillen equivalence.
\end{theorem}

\begin{proof}
We first show that $F$ reflects weak equivalences between
cofibrant objects.  Let $f \colon X \rightarrow Y$ be a map of
cofibrant Segal precategories such that $Ff \colon FX \rightarrow
FY$ is a weak equivalence of simplicial categories. (Since F
preserves cofibrations, both $FX$ and $FY$ are again cofibrant.)
Then consider the following diagram:
\[ \xymatrix{FX \ar[r] \ar[d]_\simeq & L_fFX \ar[d] & L_fX
\ar[l] \ar[d] \\
FY \ar[r] & L_fFY & L_fY. \ar[l]} \] By assumption, the leftmost
vertical arrow is a DK-equivalence.  The horizontal arrows of the
left-hand square are also DK-equivalences by definition. Since $X$
and $Y$ are cofibrant, Lemma \ref{localize} shows that the
horizontal arrows of the right-hand square are DK-equivalences.
The commutativity of the whole diagram shows that the map $L_fFX
\rightarrow L_fFY$ is a DK-equivalence and then that the map $L_fX
\rightarrow L_fY$ is also. Therefore, $F$ reflects weak
equivalences between cofibrant objects.

Now, we will show that given any fibrant simplicial category $Y$,
the map
\[ F((RY)^c) \rightarrow Y \] is a DK-equivalence. Consider a
fibrant simplicial category $Y$ and apply the functor $R$ to
obtain a Segal category which is levelwise fibrant and therefore
fibrant in $\Secat_f$.  Its cofibrant replacement will be
DK-equivalent to it in $\Secat_f$. Then, by the above argument,
strictly localizing this object will again yield a DK-equivalent
simplicial category.
\end{proof}

\section{Proofs of Lemma \ref{xnyn} and Proposition \ref{FibWE}}

In this section, we give a proof of two results stated in section
5.  We begin with a lemma which we will use in the proof of Lemma
\ref{xnyn}.

\begin{lemma} \label{comp}
Suppose that $f \colon X \rightarrow Y$ is a map of Segal
precategories with the right lifting property with respect to the
maps in $I_c$. Then
\begin{enumerate}
\item The map $f_0 \colon X_0 \rightarrow Y_0$ is surjective, and

\item The map $X_n(v_0, \ldots ,v_n) \rightarrow Y_n(fv_0, \ldots
,fv_n)$ is a weak equivalence of simplicial sets for all $n \geq
1$ and $(v_0, \ldots ,v_n) \in X_0^{n+1}$.
\end{enumerate}
\end{lemma}

\begin{proof}
Since $f \colon X \rightarrow Y$ has the right lifting property
with respect to the maps in $I_c$, it has the right lifting
property with respect to all cofibrations.  In particular, it has
the right lifting property with respect to the maps in the set
$I_f$. Therefore we can apply Lemma \ref{PQ} and the result
follows.
\end{proof}

\begin{proof} [Proof of Lemma \ref{xnyn}]
To prove Lemma \ref{xnyn}, we consider a given map $f \colon X
\rightarrow Y$ with the right lifting property with respect to the
maps in $I_c$. It follows from Lemma \ref{comp} that the map $X_0
\rightarrow Y_0$ is surjective and such that for all $n \geq 1$
and $(v_0, \ldots v_n) \in X_0^{n+1}$ the map
\[ X_n(v_0, \ldots ,v_n) \rightarrow Y_n(fv_0, \ldots ,fv_n) \]
is a weak equivalence of simplicial sets.

We must prove that $\map_{L_cX}(x,y) \rightarrow
\map_{L_cY}(fx,fy)$ is a weak equivalence of simplicial sets.
Given that fact, combining it with the surjectivity of the map
$X_0 \rightarrow Y_0$ implies that $\Ho(L_cX) \rightarrow
\Ho(L_cY)$ is an equivalence of categories.

We construct a factorization $X \rightarrow \Phi Y \rightarrow Y$
such that $(\Phi Y)_0 = X_0$ and the map
\[ (\Phi Y)_n (v_0, \ldots ,v_n) \rightarrow Y_n(fv_0, \ldots ,fv_n)
\] is an isomorphism of simplicial sets
for all $(v_0, \ldots ,v_n) \in (\Phi Y)_0^{n+1}$.  We begin by
defining the object $\Phi Y$ as the pullback of the diagram
\[ \xymatrix{\Phi Y \ar[r] \ar[d] & Y \ar[d] \\
\cosk_0(X_0) \ar[r] & \cosk_0(Y_0).} \]  Note in particular that
$(\Phi Y)_0=X_0$.

Now, notice that for each $n \geq 1$ and $(v_0, \ldots ,v_n) \in
(\Phi Y)_0^{n+1}$ the map
\[ (\Phi Y)_n(v_0, \ldots ,v_n) \rightarrow Y_n(fv_0, \ldots ,fv_n) \]
is an isomorphism of simplicial sets. Since each
\[ X_n(v_0, \ldots ,v_n) \rightarrow Y_n(fv_0, \ldots ,fv_n) \]
is a weak equivalence, we can apply model category axiom MC2 to
simplicial sets to see that the map
\[ X_n(v_0, \ldots ,v_n) \rightarrow (\Phi Y)_n(v_0, \ldots ,v_n) \]
is a weak equivalence for each $n \geq 1$ and $(v_0, \ldots v_n)$
also.

Thus we have shown that if $X \rightarrow Y$ has the right lifting
property with respect to the maps in $I_c$, then each map
$X_n(v_0, \ldots ,v_n) \rightarrow (\Phi Y)_n(v_0, \ldots ,v_n)$
is a weak equivalence of simplicial sets for $n \geq 1$ and $(v_0,
\ldots ,v_n) \in X_0^{n+1}$.  Since $X_0 = (\Phi Y)_0$, the map $X
\rightarrow \Phi Y$ is actually a Reedy weak equivalence and
therefore also a DK-equivalence. To prove Lemma \ref{xnyn}, it
remains to show that the map $\Phi Y \rightarrow Y$ is a
DK-equivalence, implying that the map $X \rightarrow Y$ is also.
We will prove this fact by induction on the skeleta of $Y$.

We will denote by $\sk_nY$ the $n$-skeleton of $Y$, as defined
above in the paragraph below Proposition \ref{inj}. We seek to
prove that the map
\[ \Phi (\sk_n Y) \rightarrow \sk_nY \] is a DK-equivalence
for all $n \geq 0$.

We first consider the case where $n=0$.  In this case, $\sk_0(\Phi
Y)$ and $\sk_0Y$ are already Segal categories. They can be
observed to be DK-equivalent as follows.  In the case of $\sk_0Y$,
given any pair of elements $(x,y) \in (\sk_0Y)_0 \times
(\sk_0Y)_0$, the mapping space $\map_{\sk_0Y}(x,y)$ is the
homotopy fiber of the map
\[ (\sk_0Y)_1 = (\sk_0Y)_0 \times (\sk_0Y)_0 \rightarrow
(\sk_0Y)_0 \times (\sk_0Y)_0 \] over $(x,y)$.  If $x=y$, this
fiber is just the point $(x,y)$, since in this case this map is
the identity. If $x \neq y$, then the fiber is empty.  For $(a,b)
\in (\sk_0 \Phi Y)_0 \times (\sk_0 \Phi Y)_0$, the fiber of the
analogous map over $(a,b)$ is equivalent to $(a,b)$ if $a$ and $b$
map to the same point $x$ in $Y_0$.  Otherwise the fiber is empty.
The definition of $\Phi Y$ and the map $\Phi Y \rightarrow Y$ then
show that the two are DK-equivalent.

We now assume that the map $\Phi (\sk_{n-1}Y) \rightarrow
\sk_{n-1}Y$ is a DK-equivalence and seek to show that the map
\[ \Phi (\sk_nY) \rightarrow \sk_nY \] is also for
$n \geq 2$.  Notice that $\sk_nY$ is obtained from $\sk_{n-1}Y$
via iterations of pushouts of diagrams of the form
\begin{equation} \label{pmnqmn}
\xymatrix@1{Q_{m,n} & \ar[l] P_{m,n} \ar[r] & \sk_{n-1}Y}
\end{equation}
For simplicity, we will assume that $m=0$ and we require only one
such pushout to obtain $\sk_nY$. Notice that
$(\sk_{n-1}Y)_0=(\sk_nY)_0=Y_0$ and that the map
\[\sk_{n-1}Y \rightarrow \sk_nY \] is the identity on the discrete space in
degree zero. Therefore we use the distinct $n$-simplex $\Delta
[n]^t_{y_0, \ldots ,y_n}$ for each $(y_0, \ldots ,y_n) \in
Y_0^{n+1}$ as defined above in section \ref{secats}. Setting $\yu
= (y_0, \ldots y_n)$, we write this $n$-simplex as $\Delta
[n]^t_\yu$.

We can then apply the map $\Phi$ to diagram \ref{pmnqmn} (and its
pushout) to obtain the diagram
\begin{equation} \label{phi}
\xymatrix{\Phi \dot \Delta [n]^t_\yu \ar[r] \ar[d] & \Phi
\sk_{n-1}Y \ar[d] \\
\Phi \Delta [n]^t_\yu \ar[r] & \Phi \sk_nY.}
\end{equation}

We would like to know that we still have a pushout diagram. In
other words, we want to know that the functor $\Phi$ preserves
pushouts. To see that it does, consider the levelwise pullback
diagram defining $(\Phi Y)_n$:
\[ \xymatrix{(\Phi Y)_n \ar[r] \ar[d] & Y_n \ar[d] \\
X_0^{n+1} \ar[r] & Y_0^{n+1}.} \]  We can regard the map $f \colon
X \rightarrow Y$ as inducing a pullback functor $f^*$ from the
category of simplicial sets over $Y_0^{n+1}$ to the category of
simplicial sets over $X_0^{n+1}$.  (Recall that the category of
objects \emph{over} a simplicial set $Z$ has as objects maps of
simplicial sets $W \rightarrow Z$ and as morphisms the maps of
simplicial sets making the appropriate triangle commute.) However,
this functor between over categories can be shown to have a right
adjoint.  Therefore it is a left adjoint and hence preserves
pushouts.

We know that the maps
\[ \Phi \dot \Delta [n]^t_\yu \rightarrow
\dot \Delta [n]^t_\yu \] and
\[ \Phi (\sk_{n-1} Y) \rightarrow \sk_{n-1}Y \]
are DK-equivalences by our inductive hypothesis, since the
nondegenerate simplices in each case are concentrated in degrees
less than $n$.  Since the left-hand vertical maps of diagrams
\ref{pmnqmn} and \ref{phi} above are cofibrations, the right-hand
vertical map in diagram \ref{phi} is also a cofibration, and
therefore it remains only to show that the map $\Phi \Delta
[n]^t_\yu \rightarrow \Delta [n]^t_\yu$ is a DK-equivalence in
order to show that the pushouts of the two diagrams are weakly
equivalent.

If $n=0$, then $\Phi \Delta [0]^t_\yu \rightarrow \Delta
[0]^t_\yu$ is a DK-equivalence since everything is already local
and $\Phi \Delta [0]^t_\yu$ is just the nerve of some contractible
category.  In fact, given any $n \geq 0$ and $\yu = (y_0, \ldots
,y_n)$, if $y_i \neq y_j$ for each $0 \leq i,j \leq n$, the map
$\Phi \Delta [n]^t_\yu \rightarrow \Delta [n]^t_\yu$ is a
DK-equivalence, since $\Delta [n]^t_\yu$ is already local.

Now suppose that $n=1$ and $\yu =(y_0, y_0)$.  Consider $g:\Phi
\Delta [1]^t_\yu \rightarrow \Delta [1]^t_\yu$ and let $k$ be the
number of 0-simplices of $g^{-1}(y_0)$.  If $C_k$ denotes the
category with $k$ objects and a single isomorphism between any two
objects, then we have that
\[ \Phi \Delta [1]^t_\yu \simeq \Delta [1]^t_\yu \times \nerve
(C_k). \] Thus, it suffices to show that
\[ L_c \Phi \Delta [1]^t_\yu \simeq L_c \Delta [1]^t_\yu \times
L_c \nerve (C_k). \]

To prove this fact, first note that the fibrant objects in
$\sesp_c$ are closed under internal hom, namely that given a Segal
space $W$ and any simplicial space $Y$, there is a Segal space
$W^Y$ given by $(W^Y)_k = \Map^h(Y \times \Delta [k]^t,W)$
\cite[7.1]{rezk}. Therefore, given any Segal precategories $X$ and
$Y$ and any Segal space $W$, we can work in the category $\sesp_c$
and make the following calculation.
\[ \begin{aligned}
\Map^h(L_cX \times L_cY, W) & \simeq \Map^h (L_cX, W^{L_cY}) \\
& \simeq \Map^h(X, W^Y) \\
& \simeq \Map^h (X \times Y, W) \\
& \simeq \Map^h (L_c(X \times Y), W)
\end{aligned} \]  In other words, the map
\[ L_c(X \times Y) \rightarrow L_cX \times L_cY \] is a
DK-equivalence, and in particular the statement above for $L_c
\Phi \Delta [n]^t_\yu$ holds.

Now consider the case where $n=2$.  Then if $\yu =(y_0, y_1,
y_2)$, we have that $G(2)^t_\yu$ can be written as a pushout
\begin{equation} \label{g}
\xymatrix{G(0)^t_{y_1} \ar[r] \ar[d] & G(1)^t_{y_0,y_1} \ar[d]
\\
G(1)^t_{y_1,y_2} \ar[r] & G(2)^t_\yu.}
\end{equation}
Now consider the map $g:\Phi G(2)^t_\yu \rightarrow G(2)^t_\yu$.
We have that $g^{-1}(G(0)^t_{y_1})$ is the nerve of some
contractible category. Similarly, the map
$g^{-1}(G(1)^t_{y_0,y_1}) \rightarrow G(1)^t_{y_0,y_1}$ is a
DK-equivalence, as is the map $g^{-1}(G(1)^t_{y_1,y_2})
\rightarrow G(1)^t_{y_1,y_2}$. Since we have a pushout diagram
\[ \xymatrix{g^{-1}(G(0)^t_{y_1}) \ar[r] \ar[d] & g^{-1}(G(1)^t_{y_0,y_1})
\ar[d] \\
g^{-1}(G(1)^t_{y_1,y_2}) \ar[r] & \Phi G(2)^t_\yu} \] and the
left-hand vertical maps of this diagram and of diagram \ref{g} are
cofibrations, it follows that the map $\Phi G(2)^t_\yu \rightarrow
G(2)^t_\yu$ is a DK-equivalence. In fact, for any $n \geq 2$,
$G(n)^t_\yu$ can be obtained by iterating such pushouts.
Therefore, we have shown that the map $\Phi G(n)^t_\yu \rightarrow
G(n)^t_\yu$ is a DK-equivalence.

To see that $\Phi \Delta [n]^t_\yu \rightarrow \Delta[n]^t_\yu$ is
a DK-equivalence for any choice of $\yu$, we need a variation on
this argument. Again using a pushout construction, we will use the
fact that this map is a DK-equivalence when each $y_i$ is distinct
to show that it is also a DK-equivalence even if $y_i =y_j$ for
some $i \neq j$.  We will describe this construction for a
specific example, but it works in general. Specifically, we show
that $\Phi \Delta [2]^t_{y_0, y_1, y_0} \rightarrow \Delta
[2]^t_{y_0, y_1, y_0}$ is a DK-equivalence.

Define the Segal precategory $\widetilde Y = Y \amalg \{\widetilde
y\}$, where $\widetilde y$ is a 0-simplex not in $Y_0$, and we
regard $\{\widetilde y\}$ as a doubly constant simplicial space.
Then, using the map $g \colon \Phi Y \rightarrow Y$ and some
vertex $y_0$ of $Y$, we let $Z$ be a Segal precategory isomorphic
to $(g^{-1}y_0)$ and define $\widetilde X= X \amalg Z$. There is a
map $\widetilde X \rightarrow \widetilde Y$ such that $Z$ maps to
$\widetilde y$. We define a functor $\widetilde \Phi$ and
factorization
\[ \xymatrix@1{\widetilde X \ar[r] & \widetilde \Phi \widetilde Y
\ar[r]^g & \widetilde Y} \] just as we defined $\Phi Y$ above.
More generally, we apply $\widetilde \Phi$ to any Segal
precategory with 0-simplices those of $\widetilde Y$ to obtain a
Segal precategory with 0-simplices those of $\widetilde X$, just
as we have been doing with $\Phi$.

Now consider the objects $G(2)^t_{y_0, y_1, \widetilde y}$ and
$\Delta [2]^t_{y_0, y_1, \widetilde y}$, each with 0-simplices
those of $\widetilde Y$.  There is a natural map
\[ G(2)^t_{y_0, y_1, \widetilde y} \rightarrow G(2)^t_{y_0, y_1,
y_0} \] where $\widetilde y \mapsto y_0$, and an analogous map
\[ \Delta [2]^t_{y_0, y_1, \widetilde y} \rightarrow \Delta
[2]^t_{y_0, y_1, y_0}. \]  We have a pushout diagram
\[ \xymatrix{G(2)^t_{y_0, y_1, \widetilde y} \ar[r] \ar[d] &
G(2)^t_{y_0, y_1, y_0} \ar[d] \\
\Delta [2]^t_{y_0, y_1, \widetilde y} \ar[r] & \Delta [2]^t_{y_0,
y_1, y_0}.} \]  Since the left-hand vertical map is a cofibration,
this map is actually a homotopy pushout diagram.

Now, from above we know that the maps
\[ \widetilde \Phi G(2)^t_{y_0, y_1, \widetilde y} \rightarrow
G(2)^t_{y_0, y_1, \widetilde y} \] and
\[ \Phi G(2)^t_{y_0, y_1, y_0} \rightarrow G(2)^t_{y_0, y_1, y_0} \]
are DK-equivalences.  We also know that the map
\[ \widetilde \Phi \Delta [2]^t_{y_0, y_1, \widetilde y} \rightarrow
\Delta [2]^t_{y_0, y_1, \widetilde y} \] is a DK-equivalence since
the 0-simplices $y_0, y_1, \widetilde y$ are distinct.  We can
consider the pushout diagram
\[ \xymatrix{\widetilde g^{-1} G(2)^t_{y_0, y_1, \widetilde y} \ar[r] \ar[d] &
g^{-1} G(2)^t_{y_0, y_1, y_0} \ar[d] \\
\widetilde g^{-1} \Delta [2]^t_{y_0, y_1, \widetilde y} \ar[r] &
\Phi \Delta [2]^t_{y_0, y_1, y_0}.} \]  which is again a homotopy
pushout diagram.  It follows that the map
\[ \Phi \Delta [2]^t_{y_0, y_1, y_0} \rightarrow \Delta [2]^t_{y_0, y_1, y_0} \] is a
DK-equivalence, completing the proof.
\end{proof}

We now proceed with the other remaining proof from section 5.

\begin{proof}[Proof of Proposition \ref{FibWE}]
Suppose that $f \colon X \rightarrow Y$ is a fibration and a weak
equivalence.  First, consider the case where $f_0 \colon X_0
\rightarrow Y_0$ is an isomorphism.  Without loss of generality,
assume that $X_0=Y_0$ and factor the map $f \colon X \rightarrow
Y$ functorially in $\SSets^{\Deltaop}_c$ as the composite of a
cofibration and an acyclic fibration in such a way that the $Y'_0$
remains a discrete space:
\[ \xymatrix@1{X \text{ } \ar@{^{(}->}[r] & Y' \ar@{->>}[r]^\simeq & Y.} \]
(We can obtain a $Y'$ with discrete zero space by taking a
factorization in $\SSets^{\Deltaop}_c$ analogous to the one we
defined for $\sesp_c$ at the beginning of section 5.) Since the
map $X \rightarrow Y$ is a DK-equivalence and the map $Y'
\rightarrow Y$ is a Reedy weak equivalence and therefore a
DK-equivalence, it follows that the map $X \rightarrow Y'$ is a
DK-equivalence.  In particular, $X \rightarrow Y'$ is an acyclic
cofibration and therefore by the definition of fibration in
$\Secat_f$ the dotted arrow lift exists in the following
solid-arrow diagram:
\[ \xymatrix{X \ar[r]^{\text{id}} \ar[d] & X \ar[d] \\
Y' \ar[r] \ar@{-->}[ur] & Y.} \] Thus, $f \colon X \rightarrow Y$
is a retract of $Y' \rightarrow Y$ and therefore a Reedy acyclic
fibration. In particular, $f$ has the right lifting property with
respect to the maps in $I_c$, since they are monomorphisms and
therefore Reedy cofibrations.

Now consider the general case, where $X_0 \rightarrow Y_0$ is
surjective but not necessarily an isomorphism. Then, as in the
proof of Lemma \ref{xnyn}, define the object $\Phi Y$ and consider
the composite map $X \rightarrow \Phi Y \rightarrow Y$. since by
the first case $X \rightarrow \Phi Y$ has the right lifting
property with respect to the maps in $I_c$, it remains to show
that $\Phi Y \rightarrow Y$ has the right lifting property with
respect to the maps in $I_c$.

Let $A \rightarrow B$ be an acyclic cofibration. Then there is a
dotted arrow lift in any solid-arrow diagram of the form
\begin{equation} \label{lift}
\xymatrix{A \ar[r] \ar[d]_\simeq & X \ar[r]^= \ar[d] & X \ar[d]
\\
B \ar[r] \ar@{-->}[urr] & \Phi Y \ar[d] \ar[r] & Y \ar[d] \\
& \cosk_0X_0 \ar[r] & \cosk_0Y_0}
\end{equation}
We would like to know that this lift $B \rightarrow X$ also makes
the upper left-hand square commute.

Suppose that $A_0 = B_0 = X_0$.  In this case, a map $B
\rightarrow Y$ together with a lifting
\[ \xymatrix{& X_0 \ar[d] \\
B_0 \ar[r] \ar@{-->}[ur] & Y_0} \] completely determines a map $B
\rightarrow \Phi Y$.  Therefore, in this fixed object set case,
there is only one possible lifting $B \rightarrow X$ in diagram
\ref{lift}, and one which makes the upper left-hand square
commute.

The map $X \rightarrow \Phi Y$ is a fibration in the fixed object
model category structure $\mathcal{LSS}ets^{\Deltaop}_{\mathcal O,
f}$ where $\mathcal O = X_0$.  However, since the cofibrations in
$\mathcal{LSS}ets^{\Deltaop}_{\mathcal O, f}$ are precisely the
monomorphisms, the acyclic fibrations are Reedy acyclic
fibrations.  Therefore, the map $X \rightarrow \Phi Y$ is a Reedy
acyclic fibration and thus has the right lifting property with
respect to all monomorphisms of simplicial spaces.  In particular,
it has the right lifting property with respect to the maps in
$I_c$.

Using the construction of $\Phi Y$ and the fact that $X
\rightarrow Y$ is a fibration and a weak equivalence, we can see
that $X_0 \rightarrow Y_0$ is surjective.  In particular, the map
$\cosk_0 X_0 \rightarrow \cosk_0Y_0$ has the right lifting
property with respect to the maps in $I_c$.  Using the universal
property of pullbacks, we can see that the map $\Phi Y \rightarrow
Y$ also has the right lifting property with respect to the maps in
$I_c$.
\end{proof}

\end{document}